\documentclass[3p]{elsarticle}

\makeatletter
\def\ps@pprintTitle{%
    \let\@oddhead\@empty
    \let\@evenhead\@empty
    \let\@evenfoot\@oddfoot
    }
\makeatother

\usepackage[table]{xcolor} % color
\usepackage{tabularx}
\usepackage{colortbl}
\usepackage{amsmath}
\usepackage{amssymb}
\usepackage{amsthm}
\usepackage{mathrsfs}
\usepackage{bm}
\usepackage{graphicx}
\usepackage[final]{hyperref}
\usepackage[ruled,vlined]{algorithm2e}
\usepackage{enumitem}
\usepackage{array}
\usepackage{multirow}
\usepackage{makecell}
\usepackage[title]{appendix}
\usepackage{empheq}
\usepackage[capitalise, nameinlink, noabbrev]{cleveref}
\usepackage{cases}
\usepackage{textcomp}
\usepackage{gensymb}
\usepackage{subcaption}
\hypersetup{
	colorlinks=true,       % false: boxed links; true: colored links
	linkcolor=darkred,        % color of internal links
	citecolor=blue,        % color of links to bibliography
	filecolor=magenta,     % color of file links
	urlcolor=blue         
}
\usepackage{float}
\usepackage{multirow}

\usepackage[colorinlistoftodos,prependcaption,textsize=normalsize]{todonotes}
\usepackage{regexpatch}

\usepackage{tcolorbox} % box

\usepackage{csvsimple}

\makeatletter
\xpatchcmd{\@todo}{\setkeys{todonotes}{#1}}{\setkeys{todonotes}{inline,#1}}{}{}
\makeatother

\usepackage{etoolbox}

\newcommand{\R}{\mathbb{R}}

\newcommand{\bbR}{\mathbb{R}}

\newcommand{\Rpow}[1]{\bbR^{#1}}

\newcommand{\grad}{\nabla}

\DeclareMathOperator*{\argmin}{arg\,min}

\DeclarePairedDelimiter\floor{\lfloor}{\rfloor}

\newtheorem{theorem}{Theorem}

\newcommand{\dd}{\,\textup{d}}

\newcommand{\dS}{\,\textup{d}S}

\newcommand{\bn}{\bm{n}}

\newcommand{\bq}{\bm{q}}

\newcommand{\bx}{\bm{x}}
\newcommand{\by}{\bm{y}}

\newcommand{\calC}{{\mathcal{C}}}

\newcommand{\calF}{{\mathcal{F}}}

\newcommand{\calL}{{\mathcal{L}}}
\newcommand{\calM}{{\mathcal{M}}}
\newcommand{\calN}{{\mathcal{N}}}

\newcommand{\calQ}{{\mathcal{Q}}}
\newcommand{\calR}{{\mathcal{R}}}

\newcommand{\calU}{{\mathcal{U}}}
\newcommand{\calV}{{\mathcal{V}}}

\newcommand{\dualDot}[2]{\langle {#1}, {#2}\rangle}

\newcommand*{\horzbar}{\rule[.5ex]{2.5ex}{0.5pt}}

\newcommand{\tbcoliii}{cyan!4}

\crefformat{appendix}{#2#1#3}

\begin{document}

\begin{frontmatter}
\title{Residual-Based Error Corrector Operator to Enhance Accuracy and Reliability of Neural Operator Surrogates of Nonlinear Variational Boundary-Value Problems$^\ast$}
\author{Prashant K. Jha$^1$\\
{$^1$Oden Institute for Computational Engineering and Sciences, The University of Texas at Austin, Austin, TX 78712, USA. \\
Email address: prashant.jha@austin.utexas.edu}\\[8pt]
$^\ast$\textbf{To appear in Computer Methods in Applied Mechanics and Engineering (DOI: \url{https://doi.org/10.1016/j.cma.2023.116595})}
}

\begin{abstract}
This work focuses on developing methods for approximating the solution operators of a class of parametric partial differential equations via neural operators. Neural operators have several challenges, including the issue of generating appropriate training data, cost-accuracy trade-offs, and nontrivial hyperparameter tuning. The unpredictability of the accuracy of neural operators impacts their applications in downstream problems of inference, optimization, and control. A framework based on the linear variational problem that gives the correction to the prediction furnished by neural operators is considered based on earlier work in JCP 486 (2023) 112104. The operator, called Residual-based Error Corrector Operator or simply Corrector Operator, associated with the corrector problem is analyzed further. Numerical results involving a nonlinear reaction-diffusion model in two dimensions with PCANet-type neural operators show almost two orders of increase in the accuracy of approximations when neural operators are corrected using the correction scheme. Further, topology optimization involving a nonlinear reaction-diffusion model is considered to highlight the limitations of neural operators and the efficacy of the correction scheme. Optimizers with neural operator surrogates are seen to make significant errors (as high as 80 percent). However, the errors are much lower (below 7 percent) when neural operators are corrected. 
\end{abstract}


\end{frontmatter}

%\tableofcontents

%%%% --------------------------- %%%%
\section{Introduction}\label{s:intro}

\begin{figure}
	\centering
	\includegraphics[width=0.6\textwidth]{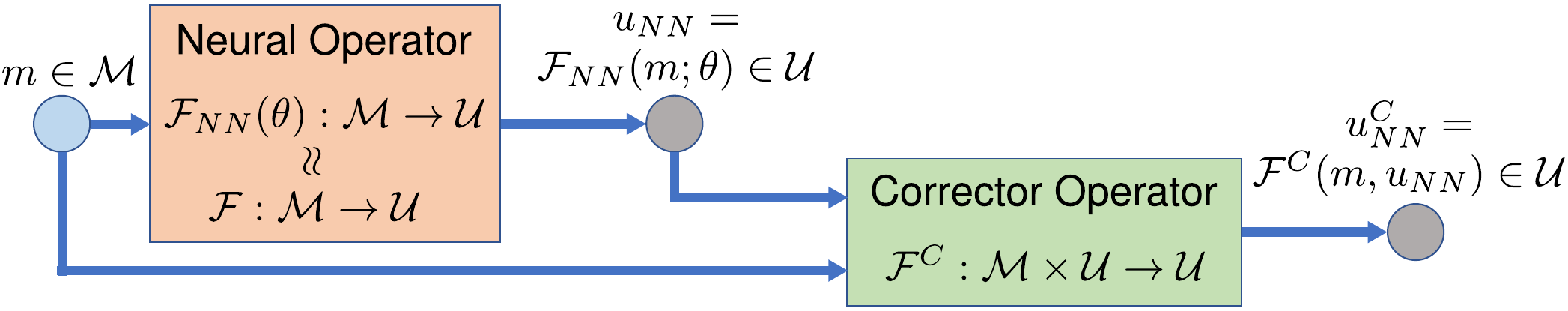}
	\caption{Schematics of the framework where the neural operator is augmented with the corrector operator developed in earlier work \cite{cao2023residual}. Here, $\mathcal{M}$ and $\mathcal{U}$ are (Banach) function spaces of model parameters and solutions of the parameterized partial differential equation, $\mathcal{F}$ the solution operator of the PDE, and $\mathcal{F}_{NN}$ a neural operator with network parameters $\theta \in \R^{d_{NN}}$. $\mathcal{F}^C$ is a correction operator defined in \cref{ss:resCorr}.}
	\label{fig:CNeuOp}
\end{figure}

This work focuses on neural operator-based surrogates constructed for a class of nonlinear parametric partial differential equations (PDEs). Specifically, neural operators that approximate the solution operator associated with the PDEs are considered. Working in a variational setting, consider a PDE (or system of PDEs) $\calR(m, u) = 0$, where $\calR$ is a residual operator, $m$ a parameter field, and $u$ a solution of PDE. Assuming that for a given $m$, there is a unique solution $u = u(m)$, an operator $\calF$ -- referred to as solution operator -- can be defined such that given $m$, $\calF(m)$ satisfies $\calR(m, \calF(m)) = 0$. For nonlinear and computationally expensive PDEs, applications in which the solution $u = \calF(m)$ is sought for large samples of parameter $m$ becomes challenging, e.g., Bayesian inference, optimization, and control under uncertainty. A neural operator $\calF_{NN}$ that maps a parameter $m$ to an approximation $\tilde{u} = \calF_{NN}(m)$ of a solution $u = u(m)$ is considered to cope with the computation cost of solving PDEs. To realize the full potential of neural operators in problems such as inference, optimization, and control possibly under uncertainties, it is essential to control the approximation errors within the required tolerance. Towards this, it is demonstrated in \cite{cao2023residual, jha2022goal} that by utilizing the underlying structure of the solution operator, an operator $\calF^C = \calF^C(\cdot, \cdot)$ can be constructed -- referred to as the Residual-based Error Corrector Operator or simply Corrector Operator -- that takes as input the parameter $m$ of the model and the neural operator prediction $\calF_{NN}(m)$ and provides a correction to the neural operator prediction; schematically shown in \cref{fig:CNeuOp}. 

In the following, the key role of parametric PDEs in various sectors is highlighted, and the neural operators are surveyed briefly. Next, the limitations of neural operators are discussed, motivating this work. The section concludes with the details of the corrector approach and notations and the layout of the paper.

\subsection{Surrogate Techniques for PDEs}

Partial Differential Equations are at the core of many engineering and scientific advancements and provide a robust and systematic means to represent physical systems or processes while encoding the fundamental
laws of mechanics such as conservation of mass, momentum, energy, and thermodynamics principles. The models of physical reality based on PDEs may include parameters (model parameters) that could change
depending on scenarios of interest or could be uncertain. As a result, one deals with a family of PDE-based models parameterized by model parameters -- so-called parametric PDEs. Examples from engineering and biomedical sectors include multiphysics modeling of complex materials  \cite{zhao2019mechanics, rahmati2023theory, darbaniyan2019designing, nandy2018monolithic, cao2022globally, lipton2019complex, jha2020kinetic, dayal2006kinetics, breitzman2018bond, lipton2018free, jha2018numerical, jha2021finite, abali2022multiphysics, dutta2021monolithic, jha2021peridynamics} and biophysical systems \cite{torbati2022coupling, hormuth2019mechanism, fritz2021analysis, fritz2021modeling, oden2016toward, oden2010general, jha2020bayesian}. 

The parametric PDE-based model is a critical component in several computationally-intensive downstream problems such as parameter estimation under uncertainty \cite{petra2014computational, hawkins2013bayesian, cao2023bayesian, karimi2023high, karimi2023hessian, bi2018bayesian, liang2023bayesian, oden2018adaptive}, topology and design optimization and optimization under uncertainty \cite{nandy2012optimization, chen2021optimal, cohen2022level, haber1996new, jog2002topology, ghattas2021learning, lipton2002optimal, lipton2002design, bendsoe1995optimal}, model selection \cite{lima2016selection, lima2017selection, luo2023optimal}, control under uncertainty \cite{chen2019taylor, alexanderian2017mean}, digital twins \cite{kapteyn2022data}, and structural health monitoring \cite{farrar2007introduction}. One common aspect of these problems is that PDE solutions are required for large samples of model parameters. If the PDEs are coupled and nonlinear, computing solutions can take significant time and resources, restricting their application in downstream problems of optimization and control. For computationally expensive PDEs, several approximation techniques are available: 
\begin{itemize}
	\item[(i)] low-fidelity approximations of the model \cite{jha2022goal, oden2002estimation}; 
	\item[(ii)] reduced-order modeling \cite{ghattas2021learning, nguyen2008efficient, qian2022reduced, geelen2023operator};
	\item[(iii)] regression techniques (e.g., polynomial chaos) to represent the quantities of interest as a function of the parameter \cite{marzouk2009dimensionality, huan2013simulation}; and 
	\item[(iv)] neural operators -- the main focus of this work -- that approximate the solution operator \cite{luo2023efficient, wu2023large, zohdi2022digital, du2023learning, goswami2020transfer, aldakheel2023efficient, wen2022u, cao2023residual, lu2021learning, goswami2020transfer, BhattacharyaHosseiniKovachkiEtAl2020,FrescaManzoni2022,KovachkiLiLiuEtAl2021,LiKovachkiAzizzadenesheliEtAl2020a,LiKovachkiAzizzadenesheliEtAl2020b,LuJinKarniadakis2019,OLearyRoseberryVillaChenEtAl2022,OLearyRoseberryDuChaudhuriEtAl2022,RaissiPerdikarisKarniadakis2019,WangWangPerdikaris2021,YuLuMengEtAl2022,Hesthaven2018, LiZhengKovachkiEtAl2021, deHoop2022cost, Hesthaven2018}.
\end{itemize}
Next, the neural operators are briefly surveyed, and their limitations are highlighted.

\subsection{Neural Operators as Surrogates of Solution Operators of PDEs}

There has been a remarkable growth in the development of neural operator-based approximations of the solution operator of parametric PDEs in recent years; for example, DeepONet \cite{WangWangPerdikaris2021, LuJinKarniadakis2019, lu2021learning}, Derivative-informed Neural Operators (DINO) \cite{OLearyRoseberryVillaChenEtAl2022}, Fourier Neural Operators (FNO) \cite{LiKovachkiAzizzadenesheliEtAl2020a, KovachkiLiLiuEtAl2021}, Graph-based Neural Operators (GNO) \cite{LiKovachkiAzizzadenesheliEtAl2020b, li2020neural}, PCA/POD-based Neural Operators (PCANet/PODNet) \cite{BhattacharyaHosseiniKovachkiEtAl2020, FrescaManzoni2022},  Physics-informed Neural Operator (PINO) \cite{LiZhengKovachkiEtAl2021}, and Wavelet Neural Operator (WNO) \cite{tripura2022wavelet}. Many of the neural operators can be characterized in a unified way; see \cite{deHoop2022cost, BhattacharyaHosseiniKovachkiEtAl2020, KovachkiLiLiuEtAl2021}. Some applications of neural operators include Bayesian inference \cite{cao2023residual}, design of materials and structures \cite{du2023learning}, digital twin \cite{zohdi2022digital}, inverse problem \cite{RaissiPerdikarisKarniadakis2019, YuLuMengEtAl2022}, multiphase flow \cite{wen2022u}, multiscale modeling \cite{aldakheel2023efficient}, optimal experimental design \cite{wu2023large}, PDE-constrained optimization \cite{luo2023efficient}, and phase-field modeling \cite{goswami2020transfer}. 

There are multiple choices of neural operators $\calF_{NN}$ for approximating the solution operator $\calF$. However, it is usually very difficult to predict and control the accuracy of the approximation, stemming from the fact that neural networks are trained to minimize an empirical error (the error is minimized in an average sense). Neural networks work well when the input parameter is in the subspace associated with the training data, and, therefore, the choice of distributions for sampling the training data becomes highly important. And, for downstream problems such as inference, optimization, and control, it is usually not possible to construct a training distribution \textit{a priori} that is representative of the parameters that may be encountered during the solution of these problems. There are techniques to remedy this issue -- the downstream problems can be solved with a crude approximation to extract crucial features of input parameters. Another option is actively updating the neural network parameters by generating new training samples. While these methods attempt to overcome the limitations of choosing appropriate training data, they may not be robust. Related to the issue of generating training data is also the trade-off between the cost and accuracy of neural operators \cite{deHoop2022cost}. 

While sufficient conditions for the existence of neural operators that are arbitrarily close to the target operator can be shown  \cite{chen1993approximations, chen1995universal}, in practice, constructing such neural operators is nontrivial. Increasing the number of training data or the complexity of neural networks may not necessarily increase the accuracy. Beyond a certain level of accuracy, it often becomes increasingly challenging to enhance accuracy, and trial-and-error approaches for tuning hyperparameters may only result in marginal gains. In this direction, adaptively increasing the complexity of neural networks and using the derivative information of the target map to create a reduced basis can help maximize the accuracy \cite{OLearyRoseberryDuChaudhuriEtAl2022}. Lastly, consider a scenario in which the training data is limited and sparse, and the neural operators employed are purely data-driven and do not explicitly enforce solving the variational problem, albeit in an approximate sense. In such scenarios, not much can be done to improve the accuracy of neural operators. 

\subsection{Approach for Enhancing Accuracy of Neural Operators}
 
Our goal towards improving the accuracy and reliability of neural operators is to utilize the underlying structure of the target map $\calF$ and build a corrector framework externally from the neural operators. A computationally inexpensive framework is sought that can take neural operator predictions and provide a correction with increased accuracy and reliability at a low cost. 

For the PDEs represented in variational form, $\calR(m, u) = 0$, so-called goal-oriented error estimates \cite{jha2022goal, oden2001goal,oden2002estimation} can provide a way forward. Following \cite{jha2022goal, oden2001goal,oden2002estimation}, given any approximation $\tilde{u}$ of the solution $u = u(m)$ of the variational problem, under certain usually reasonable assumptions, one can estimate the error $u - \tilde{u}$ by solving a linear variational problem $\delta_u \calR(m, \tilde{u}) (\tilde{e}) = -\calR(m, \tilde{u})$ for approximate error $\tilde{e}$; $\delta_u \calR(m, \tilde{u}) (\tilde{e})$ being the variational derivative of $\calR(m, \tilde{u})$ in the direction $\tilde{e}$. If $\calQ(\cdot)$ is the Quantity of Interest (QoI) functional, then it is shown in \cite{jha2022goal} that goal-oriented error, $\calQ(u) - \calQ(\tilde{u})$, can be approximated by $\delta_u \calQ(\tilde{u}) (\tilde{e})$. There are different versions of estimates available for goal-oriented error, for example, $\calQ(u) - \calQ(\tilde{u}) \approx \dualDot{\tilde{p}}{\calR(m, \tilde{u})}$, $\tilde{p}$ being the approximation of the solution of the dual problem associated with the variational form $\calR$ and the QoI functional $\calQ$; see \cite{jha2022goal, oden2001goal, oden2002estimation, van2011goal, prudhomme1999goal, prudhomme2003computable, rannacher1997feed, giles2002adjoint}

Once the error is estimated, a correction $u^C = \tilde{e} + \tilde{u}$ can be easily computed. Such an approach following \cite{jha2022goal} was developed in an earlier work \cite{cao2023residual} and it was shown in \cite{cao2023residual} that correcting predictions of trained neural operators leads to an increase in accuracy that simply can not be achieved by hyperparameter tuning or training with larger samples of data. The steps of computing estimate $\tilde{e}$ of error and the correction $u^C$ given $m$ and $\tilde{u}$ can be combined to define an operator -- referred to as the corrector operator -- $\calF^C: \calM\times \calU \to \calU$ such that $u^C = \tilde{e} + \tilde{u} = \calF^C(m, \tilde{u})$. If $\tilde{u}$ is already close to $u$, say $\tilde{u}$ is the prediction of the neural operator that has been trained to achieve a certain level of accuracy, the correction $u^C$ is expected to have two orders of more accuracy as compared to $\tilde{u}$, owing to the Newton-Kantorovich theorem; see \cref{ss:resCorr}. Further, because one only solves the linear variational problem (linear in the error estimate), the added computation cost is smaller than solving the target forward problem. 

\subsection{Contributions of This Work}
Motivated by the above observations, this work analyzes the correction scheme developed in \cite{cao2023residual} and it is shown theoretically and numerically through different examples that the correcting neural operators can have significant effects on the performance of neural operator surrogates of PDEs. It is essential to highlight here that the correction framework is external to the neural operator and uses the neural operator as a black box. To demonstrate the utility of the corrector operator, a numerical example involving a nonlinear reaction-diffusion model is considered. The accuracy of neural operators and their corrections using the corrector operator for varying input and output reduced dimensions and training sample sizes is analyzed. Moreover, topology optimization of the diffusivity parameter field in a nonlinear reaction-diffusion model is taken up to highlight the limitations of neural operators and the efficacy of the correction scheme. Particularly, three different versions of optimization problems with (1) ``true" (up to numerical discretization error) forward model, (2) neural operator surrogates of the forward model, and (3) neural operators with corrector operator are solved to compare the performance of neural operators and the improvements due to the corrector operator in the accuracy of optimizers. The results show a significant error reduction when the corrector operator is applied to neural operators. The error in the case of neural operator surrogates is as high as 80 percent while the error is seen to be below 7 percent when neural operators are corrected. 

\subsection{Notations}
Let $\mathbb{N}, \mathbb{Z}, \R$ denote the space of natural numbers, integers, and real numbers, respectively, $\R^{+}$ the space of all nonnegative real numbers. $\R^n$ denotes the $n$-dimensional Euclidean space, $x,y\in \R^n$ generic points, and $||x||$ the Euclidean norm of $x\in \R^n$. Space of $L^2$-integrable functions $f: \Omega \subset \R^{d_i} \to \R^{d_o}$ is denoted by $L^2(\Omega; \R^{d_o})$; space $H^s(\Omega; \R^{d_o})$ for functions in $L^2(\Omega; \R^{d_o})$ with generalized derivatives up to order $s$ in $L^2(\Omega; \Rpow{d_i\times_{i=1}^{s-1} d_i\, \times d_o})$. $\calL(\calU; \calV)$ denotes the space of continuous linear maps from $\calU$ to $\calV$ and $\calC^1(U; \calV)$ space of continuous and differentiable maps from $U \subset \calU$ to $\calV$.
$\calM$ and $\calU$ denote the generic Banach spaces of functions; for $u\in \calU$, $||u||_{\calU}$ denotes the norm of function $u \in \calU$. The dual of $\calU$ is the space of all linear continuous functionals on $\calU$, $L: \calU \to \R$, and is denoted by $\calU^\ast$. $\langle a, b \rangle_{\calU}$ denotes duality pairing $\calU$ and $\calU^\ast$, where $a\in \calU$ and $b\in \calU^\ast$. Given two Banach spaces $\calM$ and $\calU$, and a probability measure $\nu_M$ on $\calM$, the Bochner space of operators $\calF: \calM \to \calU$ is denoted by $L^p(\calM, \nu_{\calM}; \calU)$, for $p\in [1, \infty]$ and the norm is given by, see [Section 1.2, \cite{hytonen2016analysis}],
\begin{equation}
    ||\calF ||_{L^p(\calM, \nu_{\calM}; \calU)} = \begin{cases}
    \left(\mathbb{E}^{m \sim \nu_{\calM}} \left[ ||\calF(m) ||^p_{\calU} \right]\right)^{1/p}, &\quad p \in [1, \infty), \\
    \mathrm{ess sup}_{m \sim \nu_{\calM}} ||\calF(m)||_{\calU}, & \quad p = \infty \,.
    \end{cases}
\end{equation}
Here, $\mathbb{E}^{m \sim \nu_{\calM}} \left[ ||\calF(m) ||_{\calU} \right]$ is the expectation with respect to the probability measure $\nu_\calM$ and is defined as:
\begin{equation}
    \mathbb{E}^{m \sim \nu_{\calM}} \left[ ||\calF(m) ||_{\calU} \right] = \int_{\calM} ||\calF(m) ||_{\calU} \dd\nu_{\calM}(m) \,.
\end{equation}

\subsection{Layout of the Paper}
The paper is organized as follows: In \cref{s:pdeOp}, operators induced by PDEs in BVPs are discussed in an abstract setting, and the linear variational formulation is identified to compute the corrector in \cref{ss:resCorr}. \cref{s:neuOp} 
gives the overview of neural operators; in \cref{ss:corrScheme}, corrector framework is described. For completeness, the issues of scalability and mesh dependence with vanilla neural operators are highlighted, and the scalable approach based on singular value decomposition is discussed in \cref{ss:pcanet}. In \cref{s:numer}, numerical results are presented. The work is summarized and concluding remarks are provided in \cref{s:conc}. \cref{s:corrErrAnalysis}, \cref{s:bdResDerNonDiff}, and \cref{s:topOptDetails} include supplementary materials. 

%%%% --------------------------- %%%%
\section{Variational Boundary-Value Problems}\label{s:pdeOp}
Consider a class of parameterized variational boundary-value problems: 
\begin{equation}\label{eq:varProb}
	\begin{aligned}
		\text{Given } m\in \calM, \text{ find }u \in \calU \text{ such that }\qquad b(m, u; v) = l(v) \,,\quad\forall v\in\calU\,,
	\end{aligned}
\end{equation}
where $m$ denotes the parameter, $u$ a solution of the problem given $m$, $v$ a test function, and $\calM$ and $\calU$ are appropriate Banach function spaces associated with the parameter and solution of the problem, respectively. The semilinear form $b: \calM \times \calU \times \calU \to \bbR$ could possibly be nonlinear in the first and second arguments and linear in the last argument, and $l\in \calU^\ast$ is a continuous linear functional on $\calU$; $\calU^\ast$ being the topological dual of $\calU$. The form $b(\cdot, \cdot; \cdot)$ is assumed to characterize weak forms or variational boundary-value problems corresponding to various PDE models of physical systems or processes, with boundary conditions embedded in $\calU$ or the source term $l(\cdot)$. 

It is convenient to represent \eqref{eq:varProb} in terms of the residual operator $\calR: \calM \times \calU \to \calU^\ast$ defined via
\begin{equation}\label{eq:defResOp}
	\dualDot{v}{\calR(m, \tilde{u})} \coloneqq b(m, \tilde{u}; v) - l(v)\,, \qquad \forall v \in \calU\,,
\end{equation}
for any given $m\in \calM$ and $\tilde{u}\in \calU$. Problem \eqref{eq:varProb} now reads:
\begin{equation}\label{eq:varProbUsingResidual}
	\begin{aligned}
		\text{Given } m\in \calM, \text{ find }u \in \calU \text{ such that } \qquad\calR(m, u) = 0 \,.
	\end{aligned}
\end{equation}

Hereafter, the residual operator is assumed to be at least twice differentiable in the second argument in the variational, or, G\^ateaux sense, with the first and second derivatives, $\delta_u \calR(m, u): \calU \to \calU^\ast$ and $\delta_u^2\calR(m, u): \calU \times \calU \to \calU^\ast$ given by the linear and quadratic forms,
\begin{equation}\label{eq:defDerR}
	\begin{split}
		\delta_u \calR(m, u)(p) &:= \lim_{\epsilon \to 0} \frac{1}{\epsilon} \left[ \calR(m, u + \epsilon p) - \calR(m, u)\right]\,, \\
		\delta^2_u \calR(m, u)(p,q) &:= \lim_{\epsilon \to 0} \frac{1}{\epsilon} \left[ \delta_u \calR(m, u + \epsilon q)(p) - \delta_{u}\calR(m, u)(p)\right]\,,
	\end{split}
\end{equation}
for all $p, q\in \calU$. 

The corrector operator studied in this work is motivated by the recent work \cite{jha2022goal} on the application of so-called goal-oriented a-posteriori error estimates for the calibration of high-fidelity models using the lower fidelity approximate models. Therefore, some key aspects of goal-oriented estimates are reviewed following \cite{jha2022goal}. Following this discussion, the corrector approach is presented. 

\subsection{Goal-Oriented A-Posteriori Error Estimation}\label{ss:goal}

Consider a Quantity of Interest (QoI) $\calQ(u) \in \bbR$ to be computed using the solution $u$ of the variational problem \eqref{eq:varProbUsingResidual}. Suppose, there exists a lower fidelity model with the variational problem given by $\tilde{\calR}(\tilde{u}) = 0\, \in \calU^\ast$, $\tilde{\calR}: \calU \to \calU^\ast$ residual of lower fidelity model and $\tilde{u}$ a solution. If $\tilde{u} \approx u$, $\calQ(u)$ can be approximated by $\calQ(\tilde{u})$ with some error:
\begin{equation*}
	\calQ(u) = \calQ(\tilde{u}) + \underbrace{\calQ(u) - \calQ(\tilde{u})}_{\text{goal-oriented error}}\,.
\end{equation*}
The error $\calQ(u) - \calQ(\tilde{u})$ is referred to as the goal-oriented or modeling error as it pertains to the QoI, the goal of the analysis. Assuming that the high fidelity problem \eqref{eq:varProbUsingResidual} is computationally expensive, the goal-oriented a-posteriori error estimation provides a means to estimate goal-oriented error using $\tilde{u}$ -- a lower fidelity solution -- and some version of formula can also involve $\tilde{p}$, the dual or adjoint solution of the dual problem associated with the lower fidelity model and QoI functional $\calQ$; see \cite{jha2022goal, oden2001goal, oden2002estimation, van2011goal, prudhomme1999goal, prudhomme2003computable, rannacher1997feed, giles2002adjoint}.

Suppose $\tilde{e} = u - \tilde{u}$ is the error in the forward solution; then it is argued in [Section 2.1, \cite{jha2022goal}] that an estimate $e^C$ of $\tilde{e}$ can be computed by solving the following linear variational problem:
\begin{equation}\label{eq:approxErr}
	\delta_u \calR(m, \tilde{u})(e^C) = -\calR(m, \tilde{u})
\end{equation}
and, using $e^C$, the following estimate of the goal-oriented error can be constructed,
\begin{equation*}
	\calQ(u) - \calQ(\tilde{u}) = \delta_u \calQ(\tilde{u})(e^C) + r(u, \tilde{u}, e^C)\,,
\end{equation*}
where $r$ collects the remainder terms \cite{jha2022goal}. The equation for $e^C$ is based on the Taylor series expansion of the residual operator, as follows, see [Theorem 1.8, \cite{chow2012methods}]:
\begin{equation*}
	\calR(m, u) = \calR(m, \tilde{u} + \tilde{e}) = \calR(m, \tilde{u}) + \delta_u \calR(m, \tilde{u})(\tilde{e}) + \int_0^1 (1-s) \delta_u^2 \calR(m, \tilde{u} + s \tilde{e})(\tilde{e}, \tilde{e}) \dd s\,.
\end{equation*}
Noting that $u$ solves $\calR(m, u) = 0$, it holds that, for any $\tilde{u}\in \calU$,
\begin{equation*}
	\calR(m, \tilde{u} + \tilde{e}) = \calR(m, \tilde{u}) + \delta_u \calR(m, \tilde{u})(\tilde{e}) + \int_0^1 (1-s) \delta_u^2 \calR(m, \tilde{u} + s \tilde{e})(\tilde{e}, \tilde{e}) \dd s = 0\,.
\end{equation*}
If $\delta_{u}^2 \calR(m, \tilde{u} + s\tilde{e})$ is bounded for $s\in [0, 1]$, then the leading order term in the above is $O(||\tilde{e}||_{\calU}^2)$. It follows if $||\tilde{e}||_{\calU}$ is small, then $\tilde{e}$ satisfies $\calR(m, \tilde{u}) + \delta_u \calR(m, \tilde{e})(\tilde{e}) \approx 0 \, \in \calU^\ast$. Thus, ignoring the leading order term gives \eqref{eq:approxErr} -- an equation for the estimate $e^C$ of the error $\tilde{e}$. 

\subsection{Corrector Operator Based on Residuals}\label{ss:resCorr}

In this section, the residual-based error correction scheme developed for neural operators in an earlier work \cite{cao2023residual} based on \cite{jha2022goal} is reviewed. Suppose the variational problem \eqref{eq:varProbUsingResidual} is expensive to solve (i.e., evaluation of the solution operator $\calF(m)$), and, therefore, the operator $\calF$ is approximated using an operator $\tilde{\calF}: \calM \to \calU$ that is more easily evaluated. To ascertain how good or bad the approximation $\tilde{u}$ is, a straightforward way is to look at the exact error, $\tilde{e} = \calF(m) - \tilde{\calF}(m) = u - \tilde{u}$, and its norm, $||\tilde{e}||_{\calU}$. 
However, from the previous subsection, a computationally inexpensive method is available to estimate the error $\tilde{e}$ by $e^C$, where $e^C$ solves \eqref{eq:approxErr}. 

If $e^C$ is the estimate of error $\tilde{e} = u - \tilde{u}$, then $u^C := \tilde{u} + e^C \approx \tilde{u} + u - \tilde{u} = u$, i.e., $u^C$ is another approximation of $u$. So given an approximation $\tilde{u}$ of $u$, a linear variational problem is solved to construct another approximation $u^C$ as follows:
\begin{equation}\label{eq:corrOpVar}
	\begin{split}
		\text{Given } m\in \calM \text{ and }\tilde{u} \in \calU \text{, find } u^C \text{ such that } \quad u^C = \tilde{u} + e^C = \tilde{u} - \delta_u \calR(m, \tilde{u})^{-1}\calR(m, \tilde{u}) \,,
	\end{split}    
\end{equation}
assuming $\delta_u \calR(m, \tilde{u})^{-1}$ exists. The equation above induces an operator $\calF^C: \calM \times \calU \to \calU$, referred to as the residual-based error corrector operator or simply corrector operator, defined as
\begin{equation}\label{eq:corrOp}
	\calF^C(m, \tilde{u}) = u^C = \tilde{u} - \delta_u \calR(m, \tilde{u})^{-1}\calR(m, \tilde{u})\,,
\end{equation}
for any pair $(m, \tilde{u}) \in \calM \times \calU$. 

\subsubsection{Corrector Operator Property}

The reason $\calF^C$ is referred to as the corrector operator is that under ideal conditions, given $m\in \calM$ and an approximation $\tilde{u}$ of $u = \calF(m)$, $\calF^C$ produces another approximation $u^C = \calF^C(m, \tilde{u})$ that has smaller error as compared to $\tilde{u}$, i.e.,
\begin{equation}
    ||u - u^C||_{\calU} \leq ||u - \tilde{u}||_{\calU} \,.
\end{equation}
The following theorem provides a bound on correction error $e^C$ in terms of prediction error $\tilde{e}$.

\begin{theorem}\label{thm:corrOpConvergence}
Let $\calM$ and $\calU$ be Banach spaces, and $\calR: \calM \times \calU \to \calU^\ast$ the residual functional. For a given fixed $m \in \calM$, let $u = \calF(m)$, i.e., $u$ is such that $\calR(m, u) = 0$. For any arbitrary $\tilde{u} \in \calU$, suppose $\calR$ satisfies the following 
\begin{itemize}
    \item $\delta_u \calR(m, \tilde{u}): \calU \to \calU^\ast$ is invertible, i.e., $\delta_u \calR(m, \tilde{u})^{-1}$ exists; and
    \item $\delta^2_u \calR(m, w): \calU \times \calU \to \calU^\ast$ for all $w \in \{\tilde{u} + s (u - \tilde{u}): s \in [0, 1]\}$ is bounded.
\end{itemize}
If $u^C = \calF^C(m, \tilde{u})$, $\calF^C$ being the corrector operator defined in \eqref{eq:corrOp}, and $e^C = u - u^C$ and $\tilde{e} = u - \tilde{u}$, then the following estimate holds:
\begin{equation*}
    ||e^C||_{\calU} \leq \frac{1}{2} \left[\sup_{s\in [0,1]} \vert\vert \delta_u \calR(m, \tilde{u})^{-1} \delta^2_u \calR(m, \tilde{u} + s \tilde{e}) \vert\vert_{\calL(\calU\times \calU, \calU)} \right] \, ||\tilde{e}||_{\calU}^2\,,
\end{equation*}
where $\calL(U, V)$ is the space of all continuous linear operators from $U$ to $V$ and $|| f ||_{\calL(U, V)}$ $ = \sup_{||v||_U = 1} ||f(v)||_{V}$ is the operator norm.
\qed 
\end{theorem}
\cref{thm:corrOpConvergence} is proved in \cref{s:corrErrAnalysis}. From the above, two situations may arise:
\begin{itemize}
    \item {\it Linear error reduction.} If $\tilde{u}$ is such that $||\tilde{e}||_{\calU}$ is sufficiently small so that
    \begin{equation*}
        \frac{1}{2} \left[\sup_{s\in [0,1]} \vert\vert \delta_u \calR(m, \tilde{u})^{-1} \delta^2_u \calR(m, \tilde{u} + s \tilde{e}) \vert\vert_{\calL(\calU\times \calU, \calU)} \right] \, ||\tilde{e}||_{\calU} < 1\,,
    \end{equation*}
    then it holds $||e^C||_{\calU} < ||\tilde{e}||_{\calU}$. 

    \item {\it Quadratic error reduction.} More strictly, if $\calR$ is such that, for given $(m, \tilde{u})$ and $u = \calF(m)$,
    \begin{equation*}
        \frac{1}{2} \left[\sup_{s\in [0,1]} \vert\vert \delta_u \calR(m, \tilde{u})^{-1} \delta^2_u \calR(m, \tilde{u} + s \tilde{e}) \vert\vert_{\calL(\calU\times \calU, \calU)} \right] \leq C < 1,
    \end{equation*}
    then $||e^C||_{\calU} < ||\tilde{e}||_{\calU}^2$. Two orders of accuracy may thus be gained when $||\tilde{e}||_{\calU} < 1$. 
\end{itemize}

\paragraph{Connection to the Newton's Iteration}
The corrector operator is directly related to Newton's iteration for solving $\calR(m, u) = 0$. Let $u_0 \in \calU$ be the initial guess, then the Newton step to solve $\calR(m, u) = 0$, for a fixed $m \in \calM$, is given by
\begin{equation}
    u_{k} = u_{k-1} - \delta_u\calR(m, u_{k-1})^{-1}\calR(m, u_{k-1}) = \calF^C(m, u_{k-1}), \qquad \quad \forall k \geq 1\,.
\end{equation}
Thus, the solution $u = \calF(m)$ is a fixed-point of $\calF^C(m, \cdot)$; i.e., if $u = \calF^C(m, u)$ then $u = \calF(m)$ or, equivalently, $\calR(m, u) = 0$. Since $\calF^C$ is the operator characterized in the Newton step, the convergence of iterations, $\{u_k\}_k$, to $u$ and the reduction in error, $||u - u_{k}||_{\calU}$, as $k$ increases are ascertained by the Newton-Kantorovich theorem \cite{ortega1968newton, ciarlet2012newton}. In particular, in ideal conditions when the initial guess is sufficiently close to $u$, Newton's iterations are expected to converge at a quadratic rate, and, thus, the error reduces by two orders in every iteration. One of the versions of the Newton-Kantorovich theorem from \cite{ciarlet2012newton} is produced below.
\begin{theorem}\label{thm:newtonKantorovich}
    For a fixed $m\in \calM$, let $D(m)$ be an nonempty open set in $\calU$, $\calU$ being Banach space, and $u_0 \in D(m)$. Let $\calR(m, \cdot): \calU \to \calU^\ast$ be such that $\calR(m, \cdot) \in \calC^1(D(m); \calU^\ast)$ and $\delta_u \calR(m, u_0) \in \calL(\calU; \calU^\ast)$ is bijective, i.e., $\delta_u \calR(m, u_0)^{-1} \in \calL(U^\ast; \calU)$. Further, suppose that there exists a constant $r>0$ such that
    \begin{itemize}
        \item $\overline{B(u_0; r)} \subset D(m)$, $B(u_0; r)$ being an open ball of radius $r$ centered at $u_0$;
        \item $||\delta_u \calR(m, u_0)^{-1} \calR(m, u_0)||_{\calU} \leq \dfrac{r}{2}$; 
        \item $||\delta_u \calR(m, u_0)^{-1} \, \left(\delta_u \calR(m, \tilde{u}) - \delta_u \calR(m, \hat{u})\right)||_{\calL(\calU; \calU)} \leq \dfrac{||\tilde{u} - \hat{u}||_{\calU}}{r}$, for all $\tilde{u}, \hat{u} \in B(u_0; r)$.
    \end{itemize}
    Then, $\delta_u \calR(m, \tilde{u}) \in \calL(\calU; \calU^\ast)$ is bijective and $\delta_u \calR(m, \tilde{u})^{-1} \in \calL(\calU^\ast; \calU)$ at each $\tilde{u} \in B(u_0; r)$. The sequence $(u_k)_{k=0}^\infty$ defined by
    \begin{equation*}
        u_{k} = u_{k-1} - \delta_u \calR(m, u_{k-1})^{-1} \calR(m, u_{k-1}) = \calF^C(m, u_{k-1}), \qquad \forall k \geq 1\,,
    \end{equation*}
    is such that $u_k \in B(u_0; r)$ for all $k\geq 0$, and $u_k \to u$, where $u\in \overline{B(u_0; r)}$ is a zero of $\calR(m, \cdot)$, i.e., $\calR(m, u) = 0$. Further, for each $k \geq 0$,
    \begin{equation*}
        ||u - u_k||_{\calU} \leq \frac{r}{2^k}\,,
    \end{equation*}
    and the point $u\in \overline{B(u_0; r)}$ is the only zero of $\calR(m, \cdot)$ in $\overline{B(u_0; r)}$. 
\end{theorem}
For proof, see [Theorem 5, \cite{ciarlet2012newton}]. 

\subsection{Example of a Nonlinear Reaction-Diffusion Equation}\label{ss:nonDiff}

To put the notations and ideas discussed so far into a context, a forward problem involving a nonlinear reaction-diffusion model with homogeneous Dirichlet boundary condition is considered. Suppose $\Omega \subset \bbR^d$, $d = 1, 2, 3$, denotes the open, bounded, and smooth domain, $u = u(\bx)$, $\bx \in \Omega$, is the temperature field with $u \in \calU := H^1_0(\Omega) = \{v \in H^1(\Omega): v = 0 \text{ on }\partial \Omega\}$ governed by the reaction-diffusion model,
\begin{equation*}
	\begin{aligned}
		- \nabla \cdot (\kappa_0 m(\bx) \nabla u(\bx)) + \alpha u(\bx)^3 &= f(\bx), && \qquad\qquad\bx \in \Omega\,; \\
		u(\bx) &= 0, && \qquad\qquad \bx \in \partial \Omega \,,
	\end{aligned}
\end{equation*}
where $\kappa_0, \alpha > 0$ are fixed constants, $m \in \calM := \{v \in L^2(\Omega) \cap L^\infty(\Omega): v \geq m_{\text{lw}} \}$ is the diffusivity field, and $f \in L^2(\Omega)$ external heat source that is fixed and given. The associated variational problem reads:
\begin{equation*}
		\text{Given } m\in \calM, \text{ find }u \in \calU \text{ such that }
		\quad\underbrace{\int_{\Omega} \left\{ \kappa_0 m \nabla u \cdot \nabla v + \alpha u^3 v \right\}\dd x}_{=:\, b(m, u; v)} = \underbrace{\int_{\Omega} fv \dd x}_{=:\, l(v)}, \qquad \forall v\in \calU \,.
\end{equation*}
The corresponding residual functional $\calR: \calM \times \calU \to \calU^\ast$ is defined through action on $v\in \calU$, for $(m, u) \in \calM\times\calU$, as follows
\begin{equation*}
	\langle v, \calR(m, u)\rangle_{\calU} := b(m, u; v) - l(v) = \int_{\Omega} \left\{ \kappa_0 m \nabla u \cdot \nabla v + \alpha u^3 v \right\} \dd x - \int_{\Omega} fv \dd x \,.
\end{equation*}
The first and second derivatives, in this case, have the form:
\begin{equation}\label{eq:diffExampleResDer}
	\begin{split}
		\langle v, \delta_u\calR(m, u)(p) \rangle_{\calU} &= \int_{\Omega} \left\{ \kappa_0 m \nabla p \cdot \nabla v + 3 \alpha u^2 p v \right\} \dd x \,, \\
		\langle v, \delta^2_u\calR(m, u)(p,q) \rangle_{\calU} &= \int_{\Omega} \left\{ 6 \alpha u pq v \right\} \dd x \,.
	\end{split}
\end{equation}

The existence of solutions of the above variational problem for $d\geq 3$ can be established by following the Theorem 1.6.6 in \cite{badiale2010semilinear} and the underlying arguments. 

\subsubsection{Constants in the Corrector Operator Bound for the Nonlinear Reaction-Diffusion Equation}

The theorem below provides the bound on the first two derivatives of $\calR$ and shows that the inverse operator $\delta_u \calR(m, \tilde{u})^{-1}: \calU^\ast \to \calU$ can also be bounded. 
\begin{theorem}\label{thm:bdResDerNonDiff}
    For any $m \in \calM$ and $\tilde{u} \in \calU$, the following holds
    \begin{itemize}
        \item[(i)] Upper and lower bound on the norm of $\delta_u \calR(m, \tilde{u})$:
        \begin{equation*}
            \hat{C}_{\delta R} \,||v||_\calU \leq || \delta_u \calR(m, \tilde{u})(v) ||_{\calU^\ast} \leq \bar{C}_{\delta R}(\tilde{u}) \, ||v||_\calU\,,
        \end{equation*}
        where $\hat{C}_{\delta R}$ and $\bar{C}_{\delta R} = \bar{C}_{\delta R}(\tilde{u})$ are constants given by
        \begin{equation*}
            \hat{C}_{\delta R} := \frac{\kappa_0 m_{\text{lw}}}{2} \min\{1, C_P^{-2}\}, \qquad \bar{C}_{\delta R}(\tilde{u}) := \kappa_0 ||m||_{L^\infty(\Omega)} + 3 \alpha C_S^4 ||\tilde{u}||_{\calU}^2\,,
        \end{equation*}
        and $C_P, C_S$ are constants from the Poincar\'e inequality and the Sobolev embedding. Taking the operator norm of $\delta_u \calR(m, \tilde{u})(v)$, it holds
        \begin{equation*}
            \hat{C}_{\delta R} \leq || \delta_u \calR(m, \tilde{u}) ||_{\calL(\calU; \calU^\ast)} \leq \bar{C}_{\delta R}(\tilde{u})\,.
        \end{equation*}
        \item[(ii)] Upper bound on the norm of the inverse of $\delta_u \calR(m, \tilde{u})$:
        \begin{equation*}
            || \delta_u \calR(m, \tilde{u})^{-1} ||_{\calL(\calU^\ast; \calU)} \leq \frac{1}{\hat{C}_{\delta R}} \,.
        \end{equation*}
        \item[(iii)] Upper bound on the norm of $\delta_u^2 \calR(m, \tilde{u})$:
        \begin{equation*}
            || \delta_u^2 \calR(m, \tilde{u})(p, q)||_{\calU^\ast} \leq \bar{C}_{\delta^2 R}(\tilde{u}) \, ||p||_{\calU}\, ||q||_{\calU}, \qquad \bar{C}_{\delta^2 R}(\tilde{u}) := 6 \alpha C_S^4 ||\tilde{u}||_{\calU} \,,
        \end{equation*}
        and, therefore,
        \begin{equation*}
            || \delta_u^2 \calR(m, \tilde{u})||_{\calL(\calU \times \calU; \calU^\ast)}  \leq \bar{C}_{\delta^2 R}(\tilde{u}) \,.
        \end{equation*}
    \end{itemize}
\end{theorem}

The theorem above is proved in \ref{s:bdResDerNonDiff}. 

Using \cref{thm:bdResDerNonDiff}, the norm of the correction error, $||e^C||_{\calU}$, in terms of the norm of the prediction error $||\tilde{e}||_{\calU}$ can be estimated for the example of the nonlinear reaction-diffusion model. 
Let $\tilde{u} \in \calU$ be an approximation (prediction) of the solution $u$ of $\calR(m, u) = 0$ for a given $m \in \calM$. Let $u^C = \calF^C(m, \tilde{u})$ be the correction, where $\calF^C$ is the corrector operator defined in \eqref{eq:corrOp}. Then, combining \cref{thm:corrOpConvergence} and \cref{thm:bdResDerNonDiff}, recalling that $e^C = u - u^C$ and $\tilde{e} = u - \tilde{u}$, gives
\begin{equation*}
    \begin{split}
        ||e^C||_{\calU} &\leq \frac{1}{2} \left[\sup_{s\in [0,1]} \vert\vert \delta_u \calR(m, \tilde{u})^{-1} \delta^2_u \calR(m, \tilde{u} + s \tilde{e}) \vert\vert_{\calL(\calU\times \calU, \calU)} \right] \, ||\tilde{e}||_{\calU}^2 \\
        &\leq \frac{1}{2} \left[ \vert\vert \delta_u \calR(m, \tilde{u})^{-1} \vert\vert_{\calL(\calU; \calU^\ast)} \right]\, \left[\sup_{s\in [0,1]} \vert\vert \delta^2_u \calR(m, \tilde{u} + s \tilde{e}) \vert\vert_{\calU^\ast}\right] \, ||\tilde{e}||_{\calU}^2 \\
        &\leq \frac{1}{2\hat{C}_{\delta R}} \left[\sup_{s\in [0, 1]} \bar{C}_{\delta^2 R}(\tilde{u} + s \tilde{e}) \right] \, ||\tilde{e}||_{\calU}^2 \\
        &\leq \underbrace{\frac{3C^4_S}{\hat{C}_{\delta R}}}_{=: C} \left[ ||\tilde{u}||_{\calU} + ||\tilde{e}||_{\calU} \right] \, ||\tilde{e}||_{\calU}^2\,,
    \end{split}
\end{equation*}
i.e., under conditions in \cref{thm:bdResDerNonDiff}, a constant $C \in \R^+$ exists such that
\begin{equation*}
    ||e^C||_{\calU} \leq C ||\tilde{u}||_{\calU} \, ||\tilde{e}||_{\calU}^2 + C \, ||\tilde{e}||_{\calU}^3 \,.
\end{equation*}

%%%% --------------------------- %%%%
\section{Neural Operators and Corrector Scheme}\label{s:neuOp}
Consider a case when a map $\calF: \calM \to \calU$ is induced by a variational problem $\calR(m, \calF(m)) = 0$ in $\calU^\ast$. In practice, the variational problem is solved numerically in finite-dimensional subspaces of $\calM$ and $\calU$. In an abstract setting, the discrete variational problem can be written as:
\begin{equation}\label{eq:diffWeakFE}
	\begin{aligned}
		\text{Given } m\in \calM_h, \text{ find }u \in \calU_h \text{ such that } \qquad\dualDot{v}{\calR(m, u)} = 0\,,\quad\forall v\in\calU_h\,,
	\end{aligned}
\end{equation}
where $\calM_h \subset \calM$ and $\calU_h\subset \calU$ are finite dimensional subspaces of $\calM$ and $\calU$, respectively. For example, in a finite element approximation, if $\{\phi_i\}_{i=1}^{q_m}$ are the basis functions, where $q_m = \mathrm{dim}(\calM_h)$, then $\calM_h = \mathrm{span}\left(\{\phi_i\}\right) \subset \calM$. An element $m \in \calM_h$ is expressed as the linear combinations of interpolation functions, i.e., $ m = \sum_i m_i \phi_i$, where $(m_1, m_2, ..., m_{q_m}) \in \bbR^{q_m}$ and $\bbR^{q_m}$ being the coefficient space of $\calM_h$. 
%$\calM_h$ is isomorphic to the coefficient space $\bbR^{q_m}$. 
Similarly, if $\{\psi_i\}_{i=1}^{q_u}$ are the basis functions such that $\calU_h = \mathrm{span}\left(\{\psi_i\}\right)$ and $q_u = \mathrm{dim}(\calU_h)$, then $u \in \calU_h$ has a representation $u = \sum_{i=1}^{q_u} u_i \psi_i$ with $(u_1, u_2, ..., u_{q_u}) \in \Rpow{q_u}$. 
In what follows, $m\in \calM_h$ and $u\in \calU_h$ will be interchanged with $m = (m_1, m_2, ..., m_{q_m}) \in \bbR^{q_m}$ and $u =(u_1, u_2, ..., u_{q_u}) \in \bbR^{q_u}$ when convenient while keeping in mind that given coefficient vectors $(m_i)$ and $(u_i)$, the functions $m$ and $u$ are given by $\sum m_i\phi_i$ and $\sum u_i \psi_i$, respectively. 

The numerical discretization furnishes an approximation $\calF_h: \bbR^{q_m} \to \bbR^{q_u}$ (technically, $\calF_h: \calM_h \to \calU_h$) of $\calF: \calM \to \calU$. Next, consider a family of neural operators $\calF_{h, NN}(\cdot; \theta): \Rpow{q_m} \to \Rpow{q_u}$ parameterized by $\theta \in \Theta_{NN} \subset \Rpow{d_{NN}}$, $d_{NN}$ being the number of trainable parameters in the neural network architecture. Here, $m = (m_1, ..., m_{q_m}) \in \Rpow{q_m}$ is the input and $u = (u_1, ..., u_{q_u}) \in \Rpow{q_u}$ is the output to the neural operator. Broadly, $\theta$ is chosen such that the error $\calF_{h, NN}(\cdot; \theta) - \calF_h(\cdot)$ is minimized in some sense. 

While $\calF_h$ can be applied to any $m\in \Rpow{q_m}$, for practical purposes, some probability distribution measure $\nu$ on $\calM$ ($\nu_h$ after finite dimensional approximation of $\calM$) is assumed to sample $m$ and compute $\calF_h(m)$. Given a sampling probability distribution $\nu_h$, the optimization problem to train neural operator $\calF_{h, NN}$ can be written as:
\begin{equation}\label{eq:neuOpOptProb}
	\theta_{h, NN} = \argmin_{\theta\in \Theta_{NN}} \; J_h(\theta) := \mathbb{E}^{m \sim \nu_h} \,\left[ || \calF_h(m) - \calF_{h, NN}(m; \theta) ||_{\calU_h} \right]\,.
\end{equation}
In the above, the computation of $J_h(\theta)$ is intractable due to the nature of integration. Moreover, because $\calF_h(m)$ is expensive to compute, requiring $\calF_h(m)$ for large samples should be avoided. Thus, in practice, a finite number of samples $m^i \sim \nu_h$, $i = 1, ..., N$ (assumed independent and identically distributed) is considered to approximate the cost function as follows:
\begin{equation}\label{eq:empNeuOpOptProb}
    \tilde{\theta}_{h, NN} = \argmin_{\theta\in \Theta_{NN}} \; \tilde{J}_h(\theta) := \frac{1}{N} \sum_{i=1}^N || \calF_h(m^i) - \calF_{h, NN}(m^i; \theta) ||_{\calU_h}\,.
\end{equation}
Here, $\{(m^i, u^i = \calF_h(m^i)\}_{i=1}^N$ are the training data, and, for each $i$, $(m^i, u^i) \in \Rpow{q_m} \times \Rpow{q_u}$. Assuming the above optimization problem can be solved for $\tilde{\theta}_{h, NN}$, an ``optimal" neural operator, $\tilde{\calF}_{h, NN}(\cdot)$, can be defined according to:
\begin{equation}\label{eq:optNOp}
	\tilde{\calF}_{h, NN}(\cdot) := \calF_{h, NN}(\cdot; \tilde{\theta}_{h, NN})\,.
\end{equation}

\subsection{Correcting Neural Operators using Corrector Operator}\label{ss:corrScheme}

Suppose the optimization problem \eqref{eq:empNeuOpOptProb} is solved to obtain the neural operator $\tilde{\calF}_{h, NN}$. Since the neural operator is trained to minimize the average error, see \eqref{eq:empNeuOpOptProb}, the error $||\calF_h(m) - \tilde{\calF}_{h, NN}(m)||_{\calU_h}$ for any arbitrary sample $m \sim \nu_h$ (or $m$ can be any element of $\calM_h$) can be significantly large. Depending on the application of the neural operator-based surrogate, this unpredictability of accuracy of $\tilde{\calF}_{h, NN}(m)$ can pose a serious challenge. To improve the accuracy of $\tilde{\calF}_{h, NN}(m)$ further, one direction is to fine-tune the network architecture and hyperparameters through trial and error. However, it is seen in practice that the accuracy of a fixed neural network can not be enhanced beyond a certain limit, and often fine-tuning hyperparameters is not straightforward and may give only marginal gains \cite{cao2023residual, deHoop2022cost}. 

To enhance the accuracy and reliability beyond what can be achieved by hyperparameter tuning or increasing data samples, in an earlier work \cite{cao2023residual}, neural operator predictions are corrected using the residual-based error correction; see \cref{ss:resCorr}. This approach does not seek to modify the existing neural operator architecture and has the potential to be used with neural operators with limited accuracy and trained with sparse data. 

Given $m\in \calM_h$ and corresponding neural operator prediction $u_{NN} = \tilde{\calF}_{h, NN}(m) \in \calU_h$, the correction $u^C_{NN}$ is computed as follows:
\begin{equation*}
	u^C_{NN} = {\calF}^C(m, u_{NN}) = u_{NN} - (\delta_u \calR(m, u_{NN}))^{-1} \calR(m, u_{NN})\,.
\end{equation*}
The above entails solving the following linear variational problem:
\begin{equation}\label{eq:diffWeakCorrFE}
	\begin{aligned}
		&\text{Given } m\in \calM_h \text{ and } u_{NN} = \tilde{\calF}_{h, NN}(m) \in \calU_h, \text{ find }u^C_{NN} \in \calU_h \text{ such that } \\
		&\qquad \dualDot{v}{\delta_u \calR(m, u_{NN})(u^C_{NN} - u_{NN})} = -\dualDot{v}{\calR(m, u_{NN})}, \qquad \forall v\in \calU_h\, .
	\end{aligned}
\end{equation}
If $u_{NN}$ is sufficiently close to the true solution $u = \calF_h(m)$, as will be the case for a trained neural operator, the error $||u - u^C_{NN}||_{\calU_h}$ is expected to be at least two orders smaller than the neural operator error $||u - u_{NN}||_{\calU_h}$ when conditions of the Newton Kantorovich theorem hold; see \cref{thm:corrOpConvergence} and \cref{thm:newtonKantorovich}. 
In a scenario when $\calF_h(m)$ is sought for an input $m$ far from the subspace generated by the training input data $\{m^i\}_{i=1}^N$, the neural operator prediction is expected to have a large error, as this corresponds to extrapolation. In this case, the correction $u^C_{NN}$ is hoped to keep the error small. This is demonstrated numerically for the topological optimization of input parameter $m$ in a nonlinear reaction-diffusion equation. 
For the example of topological optimization, since it is difficult to construct \textit{a priori} a probability distribution $\nu_{h}$ that includes samples representative of $m_{optim}$, $m_{optim}$ being the solution of the topological optimization problem, the neural operator is expected to make large errors during the optimization iterations. This is the case for the topological optimization example in \cref{ss:topOpt}. However, when a corrector operator is used together with a neural operator, the accuracy of optimization solutions is seen to increase significantly. 

\subsection{Scalable and Mesh-Independent Neural Operators}\label{ss:pcanet}

The neural operator $\tilde{\calF}_{h, NN}(m): \Rpow{q_m} \to \Rpow{q_u}$ described so far, see \eqref{eq:optNOp}, has two key limitations. First, it is not scalable; for fine discretizations of variational problems, $q_m$ and $q_u$ could be large which makes training neural networks difficult due to the fact that the map is now between two very high-dimensional spaces. The second problem is mesh dependence, as $\calF_{h, NN}$ is coupled to the underlying mesh used in $\calM_h$ and $\calU_h$. 

\begin{figure}
	\centering
	\includegraphics[width=0.6\textwidth]{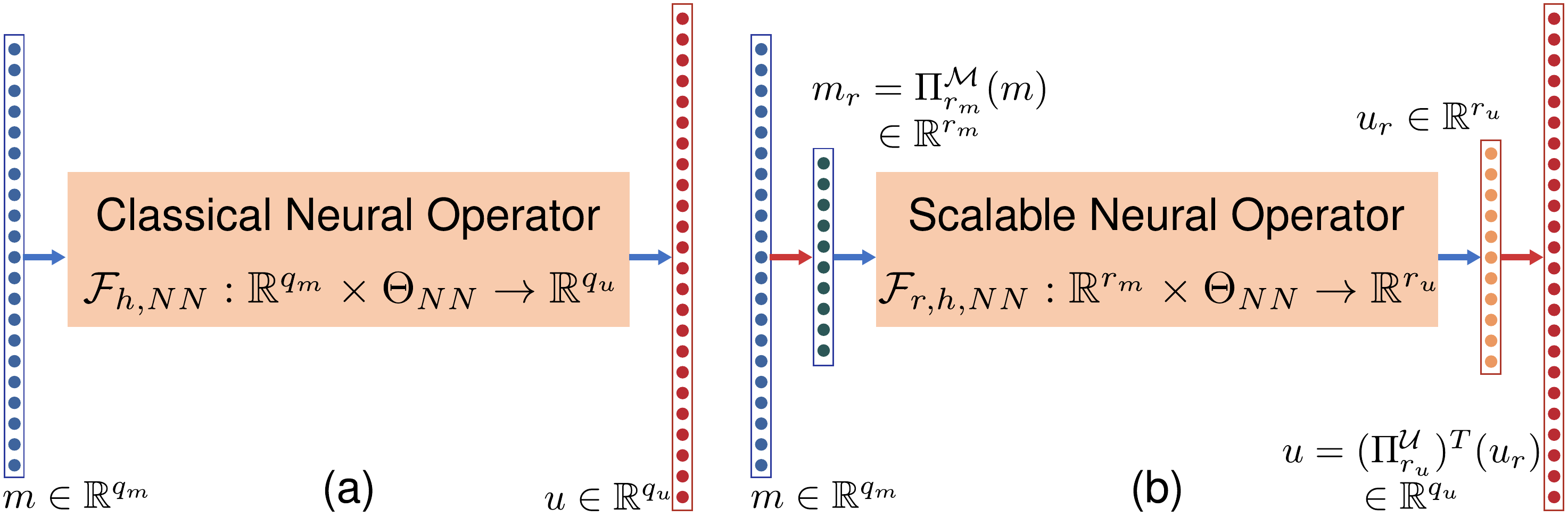}
	\caption{Schematics of the two neural operators. In (a), a classical approach is shown which maps the coefficient space of $\calM_h$ into the coefficient space of $\calU_h$. In (b), a neural operator is taken as a map between the low-dimensional subspaces of input and output spaces, and, as a result, is relatively easier to train and is independent of the mesh \cite{bhattacharya2023simulating, KovachkiLiLiuEtAl2021}. There are two additional steps in this approach: first, the given input is compressed using the projector $\Pi^{\calM}_{r_m}$, and, second, the output of the neural operator is decompressed using the transpose of the projector $\Pi^{\calU}_{r_u}$.}
	\label{fig:dimRedNOp}
\end{figure}

An approach based on dimensional reduction techniques to make the neural operators scalable and mesh-independent is discussed following \cite{BhattacharyaHosseiniKovachkiEtAl2020, OLearyRoseberryDuChaudhuriEtAl2022, OLearyRoseberryVillaChenEtAl2022, KovachkiLiLiuEtAl2021}. 
The key idea is to construct a neural network-based mapping between the low-dimensional subspaces $\Rpow{r_m}$ and $\Rpow{r_u}$ of $\Rpow{q_m}$ and $\Rpow{q_u}$, respectively; see \cref{fig:dimRedNOp}. 
To make this more precise, suppose $\Pi^{\calM}_{r_m}: \Rpow{q_m} \to \Rpow{r_m}$, $q_m = \mathrm{dim}(\calM_h)$ and $r_m \ll q_m$ is the dimension of the reduced subspace, and similarly, $\Pi^{\calU}_{r_u}: \Rpow{q_u} \to \Rpow{r_u}$, where $q_u = \mathrm{dim}(\calU_h)$ and $r_u \ll q_u$. 
Next, a parameterized neural operator $\calF_{r, h, NN}: \Rpow{r_m} \times \Theta_{NN} \to \Rpow{r_u}$ is considered with the corresponding optimization problem defined as:
\begin{equation}\label{eq:empNeuOpOptProbDiscretReduced}
	\tilde{\theta}_{r, h, NN} = \argmin_{\theta\in \Theta_{NN}} \tilde{J}_{r, h}(\theta) := \frac{1}{N} \sum_{i=1}^N || \calF_h(m_i) - (\Pi^{\calU}_{r_u})^T \left( \calF_{r, h, NN}(\Pi^{\calM}_{r_m}(m^i); \theta) \right) ||_{\calU_h}\,.
\end{equation}
The trained neural operator is then defined as
\begin{equation}\label{eq:optRedNOp}
	\tilde{\calF}_{r, h, NN}(\cdot) := \calF_{r, h, NN}(\cdot; \tilde{\theta}_{r, h, NN})\,.
\end{equation}
In the term $(\Pi^{\calU}_{r_u})^T \left( \calF_{r, h, NN}(\Pi^{\calM}_{r_m}(m^i); \theta) \right)$, firstly, the input parameter $m^i \in \Rpow{q_m}$ is projected into the reduced subspace $\Rpow{r_m}$; secondly, the reduced input vector is fed to the neural operator which returns $u_r = \calF_{r, NN}(\Pi^{\calM}_{r_m}(m^i); \theta) \in \Rpow{r_u}$; and, thirdly, $u_r$ is projected into the full space $\Rpow{q_u}$ using $(\Pi^{\calU}_{r_u})^T$. 

\subsubsection{Singular-Value Decomposition (SVD) for Projectors}\label{sss:svd}

In this work, SVD is used to construct the projectors $\Pi^{\calM}_{r_m}$ and $\Pi^{\calU}_{r_u}$ for dimensional reduction. For completeness, key aspects of SVD are reviewed in this subsection. In several works, for example, \cite{BhattacharyaHosseiniKovachkiEtAl2020}, Principal Component Analysis (PCA) provides a natural way to reduce the dimensions of discretized input and output spaces. PCA begins with the covariance or correlation matrix associated with the data matrix (see matrix $A$ below associated with the input samples for which the covariance matrix will be $\frac{1}{N-1} AA^T$) and step-by-step constructs a set of orthonormal bases that are called principal components. More precisely, the $i^\text{th}$ principal component $w^i$, such that $||w^i|| = 1$, maximizes the variance of dataset projected on $w^i$ and is orthogonal to previous principal components $w^1, w^2, ..., w^{i-1}$. In SVD, instead, one works with the data matrix directly, and the singular values and right and left singular vectors of the data matrix are sought. The singular vectors are then used as the principal bases. Following [Section 3.5, \cite{jolliffe2002principal}], the orthonormal bases computed in SVD can be shown to be related to the principal components. Therefore, the SVD can be seen as a method for performing PCA.

Suppose $\{(m^i, u^i)\}_{i=1}^N$ are the training data for the neural operator, where $m^i \in \Rpow{q_m}$ and $u^i = \calF_h(m^i)\in \Rpow{q_u}$. Further, suppose that  $\{m^i\}$ and $\{u^i\}$ are centered so that mean of $\{m^i\}$ and $\{u^i\}$, $\frac{1}{N} \sum_{i=1}^N m^i$ and $\frac{1}{N} \sum_{i=1}^N u^i$, respectively, are zero. Focusing on the input space $\Rpow{q_m}$, let $A$ denote an $q_m \times N$ matrix such that:
\begin{equation}
	A = \begin{bmatrix}
		\mid & \mid &  & \mid \\
		m^1  & m^2 & \cdots & m^N \\
		\mid & \mid & & \mid \\
	\end{bmatrix} \,.
\end{equation}
Next, consider a singular value decomposition of $A$, $A = U D V^T$, where $U$ and $V$ are column-orthonormal matrices of sizes $q_m\times q_m$ and $N \times N$, respectively, and $D$ is a $q_m\times N$ diagonal matrix. The columns of $U$ and $V$ are referred to as left and right singular vectors, respectively, while the diagonal elements of $D$, $\lambda_1 \geq \lambda_2 \geq \cdots \geq \lambda_{r} \geq 0$, are the singular values. Here,  $r = \min\{q_m, N\}$, and some $\lambda_i$ can be zero. There exists an integer $r_A \leq \min\{q_m, N\}$ such that $\lambda_{j} = 0$ for all $j > r_A$, and $r_A = \mathrm{rank}(A)$. Focusing on the matrix $U$, it has the following structure
\begin{equation}
	U = \begin{bmatrix}
		\mid & \mid &  & \mid \\
		w^1 & w^2 & \cdots & w^{q_m} \\
		\mid & \mid & & \mid \\
	\end{bmatrix} \,,
\end{equation}
where $w^i \in \Rpow{q_m}$ are orthonormal vectors, i.e., $w^i\cdot w^j = \delta_{ij}$, $\delta_{ij}$ being the Kronecker delta function. The columns of $U$, i.e., $\{w^i\}$, form a bases for $\Rpow{q_m}$.

Let $r_m > 0$ such that $r_m \leq \mathrm{rank}(A)$ is the integer of the reduced dimension $\Rpow{r_m}$ for which a projector $\Pi^{\calM}_{r_m}: \Rpow{q_m} \to \Rpow{r_m}$ is sought. Given $r_m$, a matrix $U_{r_m}$ is constructed as follows by removing the last $q_m - r_m$ columns of $U$:
\begin{equation}
	U_{r_m} = \begin{bmatrix}
		\mid & \mid & & \mid \\
		w^1 & w^2 & \cdots & w^{r_m} \\
		\mid & \mid & & \mid \\
	\end{bmatrix} \,.
\end{equation}
The matrix $U_{r_m}$ has the following notable properties:
\begin{itemize}
	\item \textit{Projection into the reduced space.} $U^T_{r_m}(m)$ projects an element $m\in \Rpow{q_m}$ into a lower dimensional subspace of $\Rpow{q_m}$, i.e., $U_{r_m}^T: \Rpow{q_m} \to \Rpow{r_m}$. To see this, consider
	\begin{equation}
		U^T_{r_m} (m) = \begin{bmatrix}
			\horzbar & (w^1)^T & \horzbar \\
			\horzbar & (w^2)^T & \horzbar \\
			& \vdots & \\
			\horzbar & (w^{r_m})^T & \horzbar \\
		\end{bmatrix} \, 
		\begin{bmatrix}
			m_1	\\
			m_2 \\
			\vdots \\
			m_{q_m} \\
		\end{bmatrix} = 
		\begin{bmatrix}
			m \cdot w^1	\\
			m \cdot w^2 \\
			\vdots \\
			m \cdot w^{r_m} \\
		\end{bmatrix} \qquad \in \Rpow{r_m} \,.
	\end{equation}
	\item \textit{Projection into the full space.} $U_{r_m}: \Rpow{r_m} \to \Rpow{q_m}$ and this is confirmed as follows: take a vector $\tilde{m} = (\tilde{m}_1, \tilde{m}_2, ..., \tilde{m}_{r_m}) \in \Rpow{r_m}$ and note
	\begin{equation}
		U_{r_m} (\tilde{m}) = \begin{bmatrix}
			\mid & \mid &  & \mid \\
			w^1 & w^2 & \cdots & w^{r_m} \\
			\mid & \mid & & \mid \\
		\end{bmatrix} \, 
		\begin{bmatrix}
			\tilde{m}_1	\\
			\tilde{m}_2 \\
			\vdots \\
			\tilde{m}_{r_m} \\
		\end{bmatrix} = 
		\sum_{i=1}^{r_m} \tilde{m}_i w^i \qquad \in \Rpow{q_m}\,,
	\end{equation}
	as $w^i \in \Rpow{q_m}$ for each $i$. 
	\item \textit{Approximation of an identity matrix.} $U_{r_m} U^T_{r_m} \approx I_{q_m}$ and $U_{r_m}^T U_{r_m} = I_{r_m}$, $I_{n}$ being the identity matrix in $\Rpow{n}$.
	\item \textit{Optimal reconstruction.} The matrix $U^T_{r_m}$ minimizes the reconstruction error over all possible projection operators of rank $r_m$. To define a reconstruction error, first note that if $m\in \Rpow{q_m}$ then $U^T_{r_m}(m) \in \Rpow{r_m}$ and $U_{r_m} (U^T_{r_m} (m)) \in \Rpow{q_m}$, i.e., $U_{r_m} (U^T_{r_m} (m))$ projects $m$ back into the same space. Thus, $U_{r_m} (U^T_{r_m} (m))$ is the reconstruction of $m$. The error, $||m - U_{r_m} U^T_{r_m} m||$, in general, is not zero. The reconstruction error for a given data matrix $A$ with columns $\{m^i\}_{i=1}^N$ is defined as the sum of the square of individual reconstruction errors:
	\begin{equation}
		e_{r_m} := \frac{1}{N} \sum_{i=1}^N ||m^i - U_{r_m}U^T_{r_m} m^i||^2 = \frac{1}{N} ||A - U_{r_m} U^T_{r_m}A||^2_{F}\,,
	\end{equation}
	where $||A||_F = \sqrt{\sum_i \sum_j A_{ij}^2}$ is the Frobenius norm of the matrix $A$. It can be shown that $U_{r_m}$ solves the following optimization problem (Eckart-Young theorem)
	\begin{equation}
		U_{r_m} = \argmin_{V \in \Rpow{q_m \times r_m}} ||A - V V^T A||_F\,.
	\end{equation}
	For proof, see [Theorem 2, \cite{chipman2020proofs}]. 
\end{itemize}
Due to the properties listed above, it makes sense to take $U^T_{r_m}$ as the projector, i.e., $\Pi^{\calM}_{r_m}:= U^T_{r_m}$. 

In a similar fashion, let $A$ is now written in terms of the output data $\{u^i\}_{i=1}^N$, i.e.,
\begin{equation}
	A = \begin{bmatrix}
		\mid & \mid & & \mid \\
		u^1  & u^2 & \cdots & u^N \\
		\mid & \mid & & \mid \\
	\end{bmatrix} \,,
\end{equation}
and $A = U D V^T$, where $U$ and $V$ are $q_u \times q_u$ and $N\times N$ orthonormal matrices, respectively, and $D$ an $q_u \times N$ diagonal matrix of singular values of $A$. 
Further, let $r_u$, such that $r_u \leq \mathrm{rank}(A) \leq \min\{q_u, N\}$, be the given dimension of the desired reduced space. The projector $\Pi^{\calU}_{r_u}: \Rpow{q_u}\to \Rpow{r_u}$ is defined using the matrix $U^T_{r_u}$, where 
\begin{equation}
	U_{r_u} = \begin{bmatrix}
		\mid & \mid & & \mid \\
		w^1 & w^2 & \cdots & w^{r_u} \\
		\mid & \mid & & \mid \\
	\end{bmatrix}
\end{equation}
is the truncation of $U$. 

%%%% --------------------------- %%%%
\section{Numerical Examples}\label{s:numer}

The example of a nonlinear reaction-diffusion equation presented in section 2.2 with slight modifications is taken up for the demonstration of the efficacy of the proposed corrector approach. The first example concerns the temperature field in a square domain with a prescribed heat source and the Dirichlet boundary condition on the bottom edge of the domain. Neural operators with varying input and output reduced dimensions and sizes of training samples are tested for accuracy. In the same tests, the corrector operator that takes neural operator prediction and model parameter as input and computes new prediction is analyzed and it is shown that the corrector operator consistently produces a new approximation of a solution of the problem with increased accuracy. The second example considers a slightly more complex geometry of a square domain with two circular voids. The forward problem now has heat flux prescribed on the outer boundary, and the temperature is fixed to zero in the inner boundaries. In this setup, the topology optimization problem on the diffusivity parameter field is posed. For this optimization problem, the accuracy of neural operators as surrogates of the forward model is examined, and it is shown that neural operators lead to high errors in optimizers. However, when the corrector operator is used in conjunction with neural operators, the accuracy increases significantly. In what follows, first, the neural network architecture and some details about the libraries used in this work are discussed. Following that, the two subsections present the key results. 

\begin{figure}[H]
	\centering
	\center
	\includegraphics[width = \textwidth]{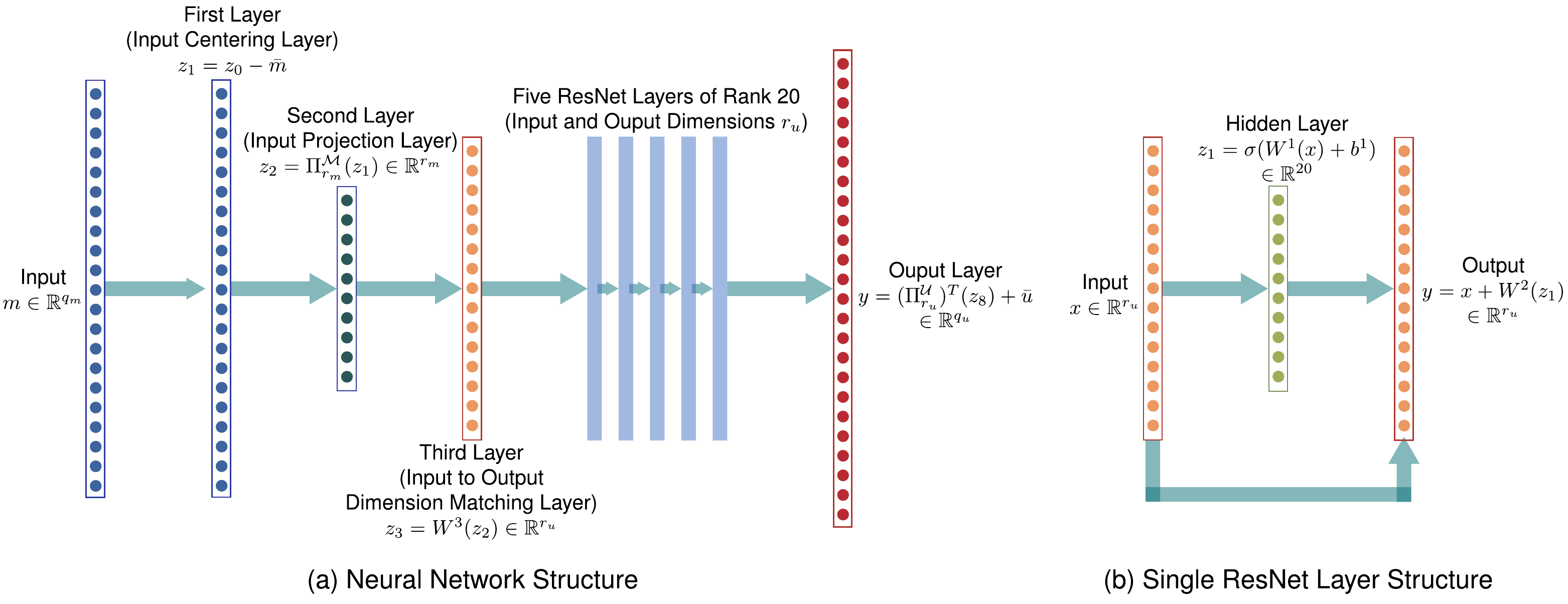}
	\caption{Neural network structure (a) and the structure of the residual network (ResNet) layer (b). In (a), data centering and projection steps are depicted. In (b), $\sigma = \sigma(z)$ is an \texttt{Softplus} ($\sigma(z) = \log(\exp(z) + 1)$) activation function. }
	\label{fig:visNet}
\end{figure}

\subsection{Neural Network Architecture and Software Details}

Neural operators in this work are based on the projectors from SVD and consist of ResNet (residual network \cite{he2016deep}) layers following \cite{OLearyRoseberryDuChaudhuriEtAl2022, cao2023residual}. Particularly, the number of residual network blocks is fixed to five with rank 20; see \cref{fig:visNet}. Let $\{(m^i, u^i)\}_{i=1}^N$ be training data with $q_u$ and $q_m$ such that $m^i \in \Rpow{q_m}$ and $u^i \in \Rpow{q_u}$, respectively, and $\bar{m} := \frac{1}{N} \sum_i m^i$ and $\bar{u} = \frac{1}{N} \sum_i u^i$. Also, let $r_m \leq q_m$ and $r_u \leq q_u$ be the dimension of reduced input and output spaces, respectively. 
In addition to the five residual network blocks, the neural network has the following four affine (an identity activation function) layers:
\begin{itemize}
	\item \textit{Input centering layer.} The first hidden layer has identity matrix as weights and $-\bar{m}$ (negative of $\bar{m}$) as bias. Given an input $m\in \Rpow{q_m}$ to this layer, the output is $m - \bar{m} \in \Rpow{q_m}$. 
	\item \textit{Input projector layer.} The second hidden layer has input projector matrix $\Pi^{\calM}_{r_m} \in \Rpow{r_m \times q_m}$ as weight, and the bias is fixed to zero. This layer projects the centered input data onto a reduced dimensional space.
	\item \textit{Input-to-output dimension matching layer.} The third layer is a \texttt{Dense} layer with $r_u$ neurons which takes $\Rpow{r_m}$ element and outputs $\Rpow{r_u}$ element. The bias is fixed to zero. 
	\item \textit{Output centering and projector layer.} The last (output) layer consists of $(\Pi^{\calU}_{r_u})^T \in \Rpow{q_u \times r_u}$ as weight and $\bar{u}$ as bias. This layer takes the output $u \in \Rpow{r_u}$ of the second last layer and projects it onto the output space $\Rpow{q_u}$ and translates by $\bar{u}$.
\end{itemize}
The weights and biases of the above layers except for the third layer for handling dimension mismatch are frozen, however, it is possible to learn the projectors (weights) by making these layers part of training \cite{OLearyRoseberryDuChaudhuriEtAl2022}. The resulting neural network with the above three layers and the five hidden layers based on ResNet is depicted in \cref{fig:visNet}. The parameters that will be varied in the numerical examples are dimensions of reduced spaces ($r_m$ and $r_u$) and the number of data ($N$). The $N$ samples of data are divided into $\floor{0.1N}$ number of validation data $N - \floor{0.1N}$ number of training data, and the testing data sample in addition to $N$ training samples is fixed to $\floor{0.25N}$. The implementation of neural networks is based on hIPPYflow\footnote{\url{https://github.com/hippylib/hippyflow}} \cite{OLearyRoseberryVillaChenEtAl2022, OLearyRoseberryDuChaudhuriEtAl2022, cao2023residual} and TensorFlow\footnote{\url{https://www.tensorflow.org/}} \cite{tensorflow2015-whitepaper}. To solve the variational problems and sample from a prior $\nu_h$,  FEniCS\footnote{\url{https://fenicsproject.org/}} \cite{alnaes2015fenics, alnaes2014unified} and hiPPYlib\footnote{\url{https://github.com/hippylib/hippylib}} \cite{villa2018hippylib} are used. 

\subsection{Accuracy Comparison for a Nonlinear Reaction-Diffusion Equation}\label{ss:prob1}

Consider a square domain $\Omega = (0, 1)^2$ with the bottom edge denoted by $\Gamma_b = [0,1]\times {0}$. The equation for the temperature field $u = u(\bx)$ over the domain $\Omega$ is taken as
\begin{equation}\label{eq:diffStrongProb1}
	\begin{aligned}
		- \nabla \cdot (e^{m(\bx)} \nabla u(\bx)) + u(\bx)^3 &= f(\bx), && \qquad \qquad\bx \in \Omega\,; \\
		u(\bx) &= 0, && \qquad \qquad \bx \in \Gamma_b \,; \\
		e^{m(\bx)} \grad u(\bx) \cdot \bn(\bx) &= 0, && \qquad \qquad \bx \in \partial \Omega - \Gamma_b\,;
	\end{aligned}
\end{equation}
where 
\begin{equation*}
	f(\bx) = e^{-4(1-x_1)^2} \sin(4 \pi x_2)^2, \qquad \qquad \bx = (x_1, x_2) \in \Omega,
\end{equation*}
is an external heat source. Let the parameter and solution function spaces be given by:
\begin{equation*}
	m \in \calM \coloneqq L^2(\Omega)\cap L^\infty(\Omega), \qquad\qquad  u \in \calU \coloneqq \{v\in H^1(\Omega): v(\bx) = 0, \, \bx \in \Gamma_b\}\,.
\end{equation*}
The variational problem associated with \eqref{eq:diffStrongProb1} reads:
\begin{equation}\label{eq:diffWeakProb1}
	\begin{aligned}
		&\text{Given } m\in \calM, \text{ find }u \in \calU \text{ such that }\\
		&\quad\dualDot{v}{\calR(m, u)}\coloneqq\int_{\Omega} e^{m(\bx)} \grad u(\bx) \cdot\grad v(\bx)\,\dd \bx + \int_\Omega u(\bx)^3 v(\bx)\,\dd\bx - \int_{\Omega} f(\bx) v(\bx)\,\dd \bx = 0\,,\quad\forall v\in\calU\,.
	\end{aligned}
\end{equation}
Expressions for $\delta_u \calR(m, u)(p)$ and $\delta_u^2 \calR(m, u)(p, q)$ can be derived following \cref{ss:nonDiff}. 

\subsubsection{Data Generation and Neural Operators}\label{sss:dataGen}
A bilinear finite element space $\calU_h$ on a quadrilateral mesh of $\Omega$ consisting of 64$\times$64 elements is considered. $\calM_h$ is identical to $\calU_h$. By replacing $\calM$ and $\calU$ with $\calM_h$ and $\calU_h$, respectively, the discrete version of the problem is obtained. 

\begin{figure}[h]
	\begin{subfigure}[t]{.5\textwidth}
		\centering
		\includegraphics[width = 0.9\textwidth]{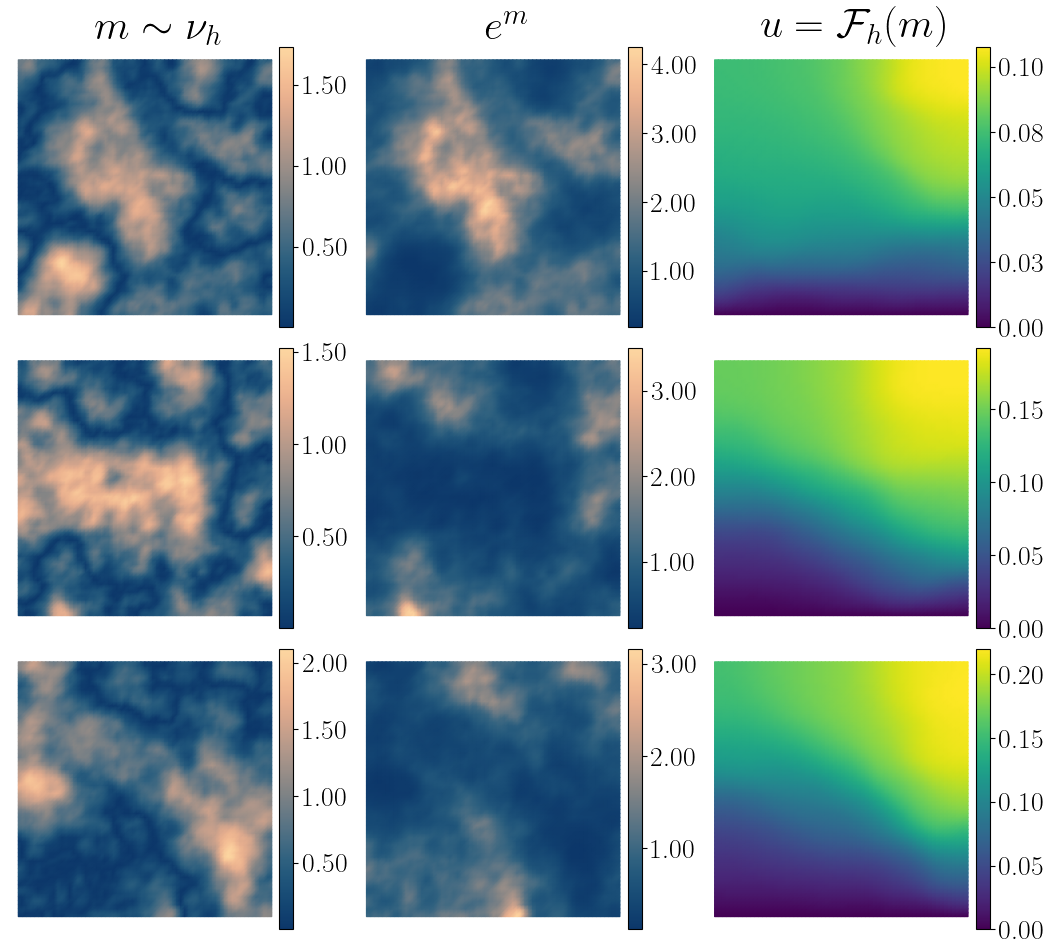}
		%\vspace{-10pt}
		\caption{}
	\end{subfigure}%
	\begin{subfigure}[t]{.5\textwidth}
		\centering
		\includegraphics[width = \textwidth]{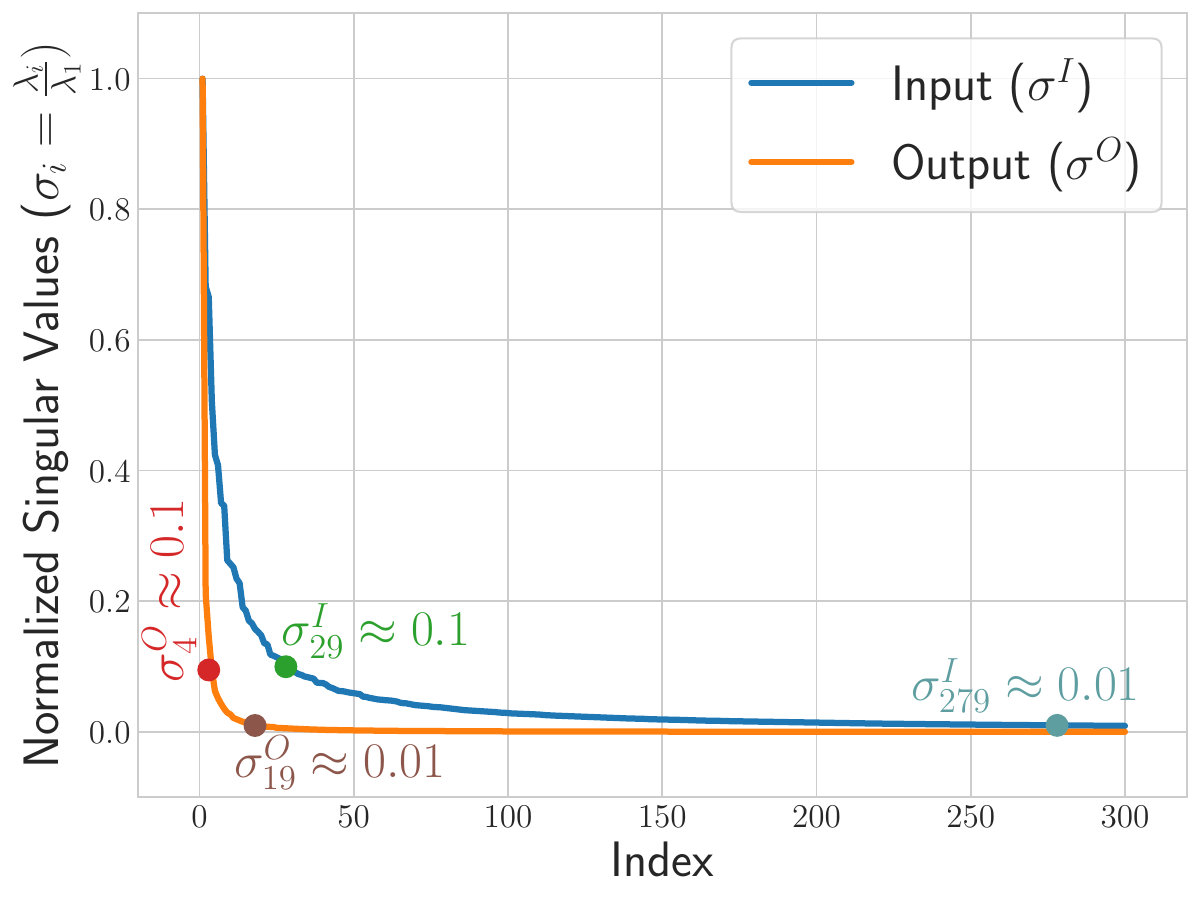}
		%\vspace{-10pt}
		\caption{}
	\end{subfigure}
	\vspace{-5pt}
	\caption{(a) Visualization of three representative samples. (b) Normalized singular values for the input and output data of sample size $N = 4096$. The indices associated with the singular values near $0.1$ and $0.01$ are also shown.}
	\label{fig:diffSamples}
\end{figure}

\paragraph{Data for Neural Operator Training}
To generate training data, probability distribution $\nu$ is assumed to be $\nu = \calN(0, \calC)$, where $\calC: \calM \times \calM \to \bbR$ is a covariance operator taking the form
\begin{equation*}\label{eq:gaussian_prior}
	\calC = \begin{cases}
		(-\gamma\nabla\cdot\nabla + \delta)^{-d} &\quad \text{ in }\Omega\,, \\
		\gamma\bn\cdot\grad + \eta &\quad \text{ on }\partial\Omega\,,
	\end{cases}
\end{equation*}
where $\gamma, \delta, \eta, d$ are hyperparameters of a covariance operator and $\bn$ unit outward normal. Covariance parameters are fixed as follows: $\gamma = 0.08, \delta = 2, \eta = 1/1.42, d = 2$. Let $\calC_h: \calM_h \times \calM_h \to \bbR$ is the covariance operator in the finite dimensional setting and $\nu_h = \calN(0, \calC_h)$. 
The set $\{(m^i, u^i)\}_{i=1}^N$, where $m^i \sim \nu_h \in \calM_h$ and $u^i = \calF_h(m^i)$ is the solution of the discretized variational problem, is the data for neural operator learning. In \cref{fig:diffSamples}(a), three representative data samples along with the singular values of input and output data from $N = 4096$ samples are depicted. 

\paragraph{Neural Operators}

To test the effect of sample size and approximations due to dimension reductions on the accuracy of neural operators, neural operators with 
\begin{equation*}
	(r_m, r_u) \in \{(50, 25), (50, 50), (100, 25), (100, 50)\} \text{ and } N \in \{256, 512, 1024, 2048, 4096\}
\end{equation*}
are trained. Here, $r_u$ is kept small compared to $r_m$ based on the relatively faster decay of output singular values; see \cref{fig:diffSamples}(b). Thus, a total of 20 neural operators are trained and tested.

\subsubsection{Comparing Neural Operator and Corrector Operator Accuracy}

\begin{figure}[h]
	\centering
	\vspace{-5pt}
	\includegraphics[width = 0.6\textwidth]{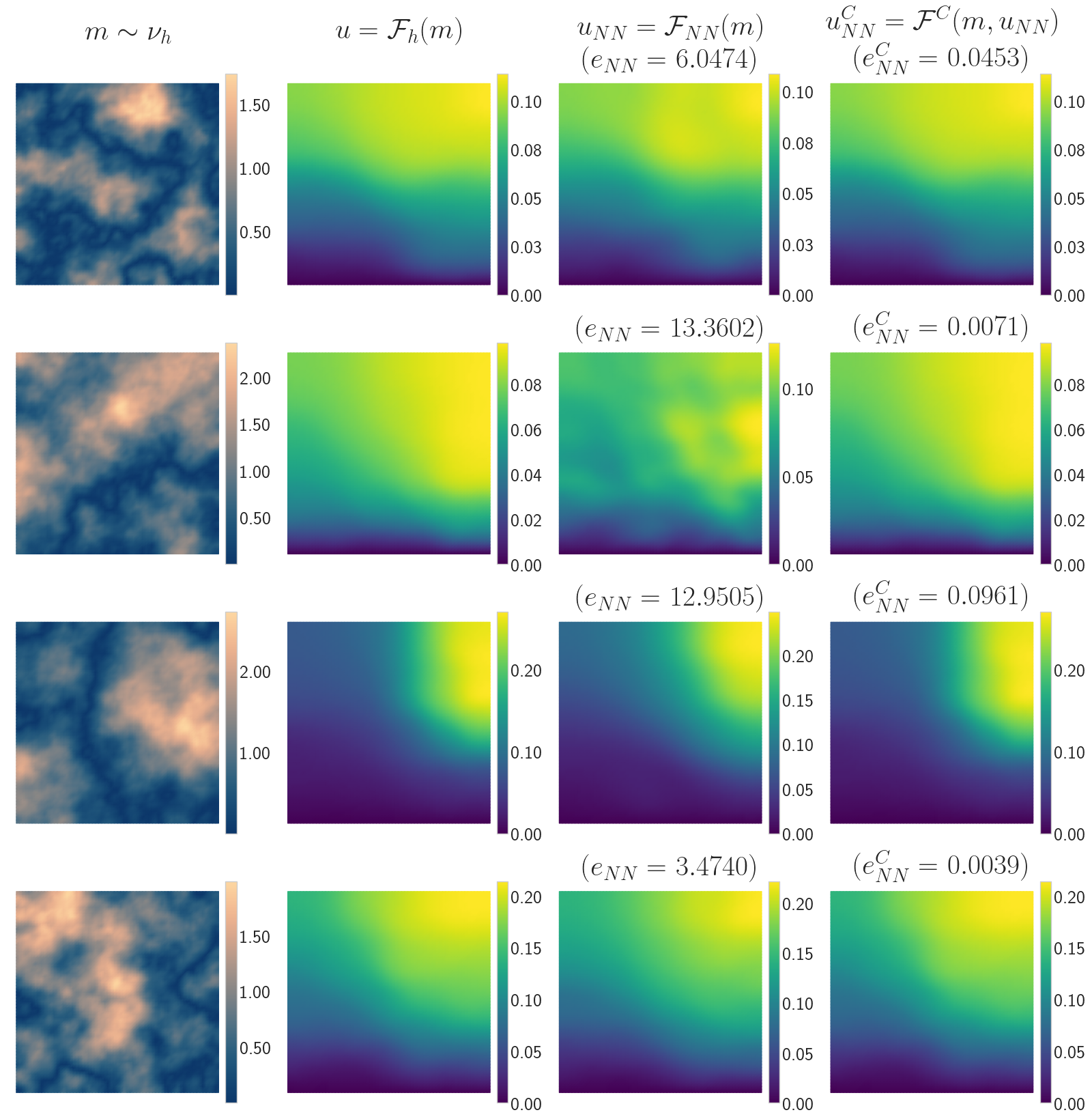}
	\vspace{-5pt}
	\caption{Comparing true solution, neural operator prediction, and the correction of neural operator prediction for networks $1, 6, 11, 16$ (see \cref{tab:compareAcc}) trained with smaller samples of data.}
	\label{fig:compareAccVis1}
\end{figure}

\begin{figure}[h]
	\centering
	\vspace{-5pt}
	\includegraphics[width = 0.6\textwidth]{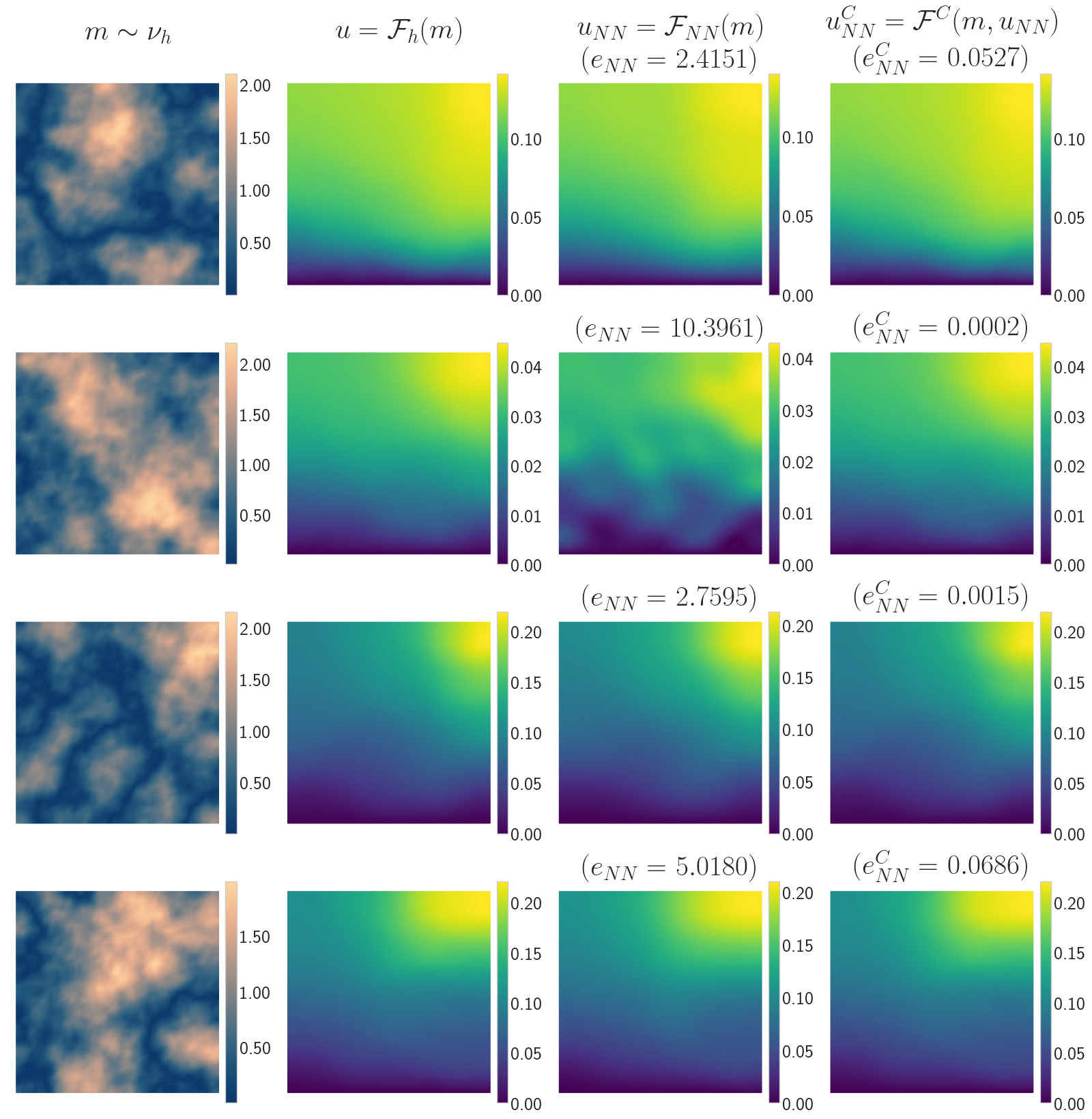}
	\vspace{-5pt}
	\caption{Comparing true solution, neural operator prediction, and the correction of neural operator prediction for networks $5, 10, 15, 20$ (see \cref{tab:compareAcc}) trained with larger samples of data.}
	\label{fig:compareAccVis2}
\end{figure}

Given a sample of input parameter $m\in \Rpow{q_m}$, suppose $u = u(m) \in \Rpow{q_u}$ is the finite element solution, $u_{NN} = u_{NN}(m)$ is the approximation furnished by neural operator, and $u^C_{NN} = u^C_{NN}(m)$ is the correction of $u_{NN}$ obtained through corrector operator. The normalized percentage error can be defined as
\begin{equation}\label{eq:errsNormalized}
	e_{NN}(m) := \frac{||u(m) - u_{NN}(m)||_{l^2}}{||u(m)||_{l^2}} \times 100\,, \qquad \qquad e^C_{NN}(m) := \frac{||u(m) - u^C_{NN}(m)||_{l^2}}{||u(m)||_{l^2}} \times 100\,,
\end{equation}
where $||a||_{l^2} = \sqrt{\sum_i a_i^2}$. In \cref{tab:compareAcc}, the statistics of errors due to neural operators and corrector is shown. Analyzing the highlighted columns in \cref{tab:compareAcc} corresponding to the mean of $e_{NN}$ and $e^C_{NN}$, the corrector approach is seen to consistently decrease the errors. In fact, for neural operators trained on small data (see rows with Numbers 1, 6, 11, and 16), the corrector does a great job of keeping the average error below $0.1$ percentage. For the neural operators trained with the smallest and largest datasets, samples of $m\sim \nu_h$ are drawn randomly, and the solutions from the true model, neural operator, and correction of the neural operator are visualized in \cref{fig:compareAccVis1} and \cref{fig:compareAccVis2}. 

{
	\begin{table}
		\centering
		{\footnotesize
			%\rowcolors{3}{\tbcoli}{\tbcolii}
			\begin{tabular}{|c|c|c|c|c|c|>{\columncolor{\tbcoliii}}c|c|c|>{\columncolor{\tbcoliii}}c|}
				\hline
				&&&& \multicolumn{3}{c|}{$e_{NN}$} & \multicolumn{3}{c|}{$e^C_{NN}$} \\ \hline 
				Number & $r_m$ & $r_u$ & $N$ & $\mathrm{min}$ & $\mathrm{max}$ & $\mathrm{mean}$ & $\mathrm{min}$ & $\mathrm{max}$ & $\mathrm{mean}$ \\
				\hline
				\rowcolor{blue!10} 
				1 & 50 & 25 & 256 & 2.45728 & 16.84489 & \cellcolor{\tbcoliii} 6.06698 & 0.00024 & 0.31488 & \cellcolor{\tbcoliii} 0.04303 \\
				2 & 50 & 25 & 512 & 1.40506 & 7.82594 & 3.53658 & 0.00016 & 0.16437 & 0.03716 \\
				3 & 50 & 25 & 1024 & 1.65162 & 6.06503 & 3.18435 & 0.00005 & 0.09237 & 0.01783 \\
				4 & 50 & 25 & 2048 & 1.40497 & 4.91186 & 2.56645 & 0.00141 & 0.20477 & 0.07224 \\
				\rowcolor{blue!10}
				5 & 50 & 25 & 4096 & 1.16622 & 9.80585 & \cellcolor{\tbcoliii} 2.43871 & 0.00001 & 0.19827 & \cellcolor{\tbcoliii} 0.03807 \\ \hline
				\rowcolor{green!10}
				6 & 50 & 50 & 256 & 3.86182 & 29.05746 & \cellcolor{\tbcoliii} 8.91953 & 0.00008 & 0.24863 & \cellcolor{\tbcoliii} 0.07290 \\
				7 & 50 & 50 & 512 & 2.89625 & 10.76244 & 5.54185 & 0.00014 & 0.19987 & 0.05248 \\
				8 & 50 & 50 & 1024 & 1.92660 & 8.37186 & 3.77333 & 0.00036 & 0.28288 & 0.06235 \\
				9 & 50 & 50 & 2048 & 1.88078 & 4.93965 & 3.07723 & 0.00011 & 0.18041 & 0.03552 \\
				\rowcolor{green!10}
				10 & 50 & 50 & 4096 & 1.74788 & 4.84569 & \cellcolor{\tbcoliii} 3.19339 & 0.00012 & 0.20722 & \cellcolor{\tbcoliii} 0.06781 \\ \hline
				\rowcolor{red!10}
				11 & 100 & 25 & 256 & 2.76545 & 11.96616 & \cellcolor{\tbcoliii} 5.23490 & 0.00937 & 0.20000 & \cellcolor{\tbcoliii} 0.05499 \\
				12 & 100 & 25 & 512 & 1.74937 & 5.60931 & 3.70606 & 0.00072 & 0.26213 & 0.07850 \\
				13 & 100 & 25 & 1024 & 1.65959 & 6.97599 & 3.37611 & 0.00019 & 0.23674 & 0.05536 \\
				14 & 100 & 25 & 2048 & 1.01940 & 5.27190 & 2.60825 & 0.00012 & 0.20808 & 0.06776 \\
				\rowcolor{red!10}
				15 & 100 & 25 & 4096 & 1.22295 & 4.00625 & \cellcolor{\tbcoliii} 2.00097 & 0.00005 & 0.27715 & \cellcolor{\tbcoliii} 0.05748 \\ \hline
				\rowcolor{orange!10}
				16 & 100 & 50 & 256 & 2.30358 & 23.87869 & \cellcolor{\tbcoliii} 5.58995 & 0.00033 & 0.23157 & \cellcolor{\tbcoliii} 0.05828 \\ 
				17 & 100 & 50 & 512 & 2.26997 & 12.24735 & 5.29943 & 0.00077 & 0.27386 & 0.05864 \\
				18 & 100 & 50 & 1024 & 1.39031 & 7.22452 & 3.02298 & 0.00008 & 0.25664 & 0.06455 \\
				19 & 100 & 50 & 2048 & 1.08747 & 6.17697 & 3.15330 & 0.00012 & 0.24090 & 0.06358 \\
				\rowcolor{orange!10}
				20 & 100 & 50 & 4096 & 1.00965 & 4.63677 & \cellcolor{\tbcoliii} 2.27597 & 0.00005 & 0.20564 & \cellcolor{\tbcoliii} 0.07560 \\
				\hline
			\end{tabular}
		}
		\caption{Comparing errors due to the neural operator approximations and corrections of neural operators for the first example (see \cref{ss:prob1}). For each neural operator and the corrector of the neural operator, the errors are computed for a total of twenty samples, and the minimum, maximum, and mean were computed from the resulting twenty $e^C_{NN}$ and $e^C_{NN}$ values. The errors are defined as in \eqref{eq:errsNormalized}. Particularly, columns corresponding to the mean of $e_{NN}$ and $e^C_{NN}$ errors are highlighted; from the results, the corrector is seen to consistently enhance the accuracy of neural operators by almost two orders. The neural operators trained with the smallest and largest data samples for different $(r_m, r_u)$ pairs are highlighted. Generally, increasing the sample size increases the accuracy of neural operators, as seen from the results.} \label{tab:compareAcc}
	\end{table}
}

%%%% --------------------------- %%%%
\subsection{Topology Optimization Involving a Nonlinear Reaction-Diffusion Equation}\label{ss:topOpt}

\begin{figure}[h]
	\center
	\includegraphics[width = 0.3\textwidth]{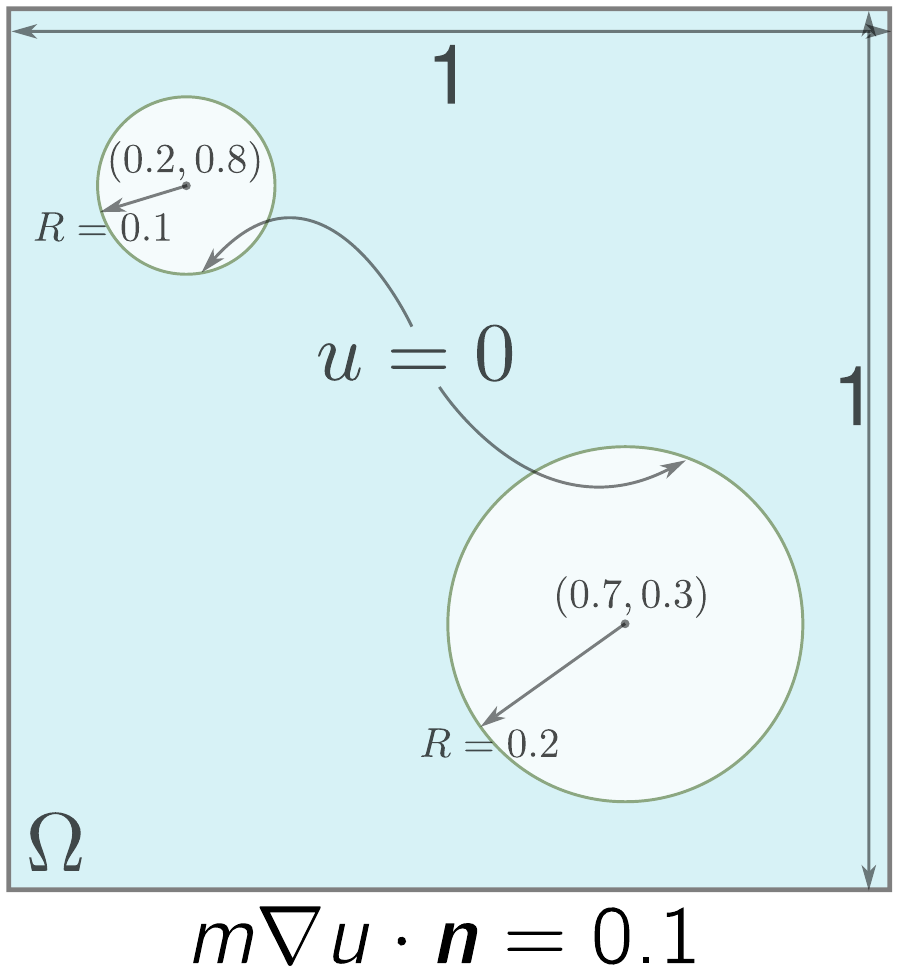}
	\caption{Setup for the second example of topological optimization of diffusivity in nonlinear reaction-diffusion equation.}
	\label{fig:diffSetup}
\end{figure}

To further test the utility of the predictor-corrector approach and highlight the limitation of neural operators in optimization problems, topological optimization of the diffusivity field in a nonlinear reaction-diffusion model is considered in this subsection. The domain $\Omega$ is a square domain with two circular voids: $\Omega = (0, 1)^2 - \overline{B(\bx_{c_1}, R_1)} - \overline{B(\bx_{c_2}, R_2)}$, where $B(\bx, R) = \{\by\in \Rpow{2}: |\by - \bx| < R \}$ denotes the ball of radius $R$ centered at $\bx$. Here, $\bx_{c_1} = (0.2, 0.8)$, $\bx_{c_2} = (0.7, 0.3)$, $R_1 = 0.1$, and $R_2 = 0.2$; see \cref{fig:diffSetup}. Let $\partial \Omega = \Gamma_{in}\cup \Gamma_{out}$, $\Gamma_{in}$ and $\Gamma_{out}$ being the inner and outer boundaries, respectively. In the inner boundary, $\Gamma_{in}$, temperature is fixed to zero, while, in the outer boundary, $\Gamma_{out}$, the heat flux $g(\bx):= 0.1$ is prescribed. Keeping the model same as in \eqref{eq:diffStrongProb1}, but now with heat source zero, $f = 0$, the strong form of the forward model reads as: 
\begin{equation}\label{eq:diffStrongProb2}
	\begin{aligned}
		\text{Given a diffusivity field } m = m(\bx) \text{, find } &\text{temperature } u \text{ such that} \\
		- \nabla \cdot ({m(\bx)} \nabla u(\bx)) + u(\bx)^3 &= 0, && \qquad \qquad\bx \in \Omega\,; \\
		u(\bx) &= 0, && \qquad \qquad \bx \in \Gamma_{in} \,; \\
		{m(\bx)} \grad u(\bx) \cdot \bn(\bx) &= 0.1 =: g(\bx), && \qquad \qquad \bx \in \Gamma_{out}\,.
	\end{aligned}
\end{equation}
As before, function spaces associated with the parameter and solution are taken as:
\begin{equation}\label{eq:fnSpacesTopo}
	m \in \calM \coloneqq L^2(\Omega) \cap L^\infty(\Omega), \qquad \qquad \quad u \in \calU \coloneqq \{v\in H^1(\Omega): v(\bx) = 0, \, \bx \in \Gamma_{in}\}\,.
\end{equation}
The variational problem corresponding to the forward problem reads:
\begin{equation}\label{eq:diffWeakProb2}
	\begin{aligned}
		&\text{Given } m\in \calM, \text{ find }u \in \calU \text{ such that }\\
		&\dualDot{v}{\calR(m, u)}\coloneqq\int_{\Omega} {m(\bx)} \grad u(\bx) \cdot\grad v(\bx)\,\dd \bx + \int_\Omega u(\bx)^3 v(\bx)\,\dd\bx - \int_{\Gamma_{out}} g \, v(\bx)\,\dS (\bx) = 0\,,\quad\forall v\in\calU\,.
	\end{aligned}
\end{equation}

\paragraph{Topological Optimization Problem}

Given a temperature $u$ satisfying the above equation, the compliance -- external working -- is defined as
\begin{equation}\label{eq:compliance}
	J(m) = \int_{\Gamma_{out}} g \, u(\bx) \dS (\bx)\,.
\end{equation}
In this example, diffusivity $m$ is optimized to minimize the compliance $J$. Let $\calM_{ad} =\{m \in \calM: 0 < m_{\text{lw}} \leq m \leq 1\} \subset \calM$ be the admissible space, $m_{\text{lw}} > 0$ being a small number suitably chosen to ensure wellposedness of the problem \eqref{eq:diffWeakProb2}. Further, let $\eta \in (0, 1]$ is the target average diffusivity, $g = 0.1$ an external heat flux on $\Gamma_{out}$, and $\calF$ a forward solution operator. The topology optimization problem reads
\begin{equation}\label{eq:optProb}
	\begin{aligned}
		\tilde{m} = &\argmin_{m \in \calM_{ad}} J(m) := \int_{\Gamma_{out}} g \, \calF(m) \dS(\bx) && \text{such that} \qquad \frac{1}{|\Omega|} \int_{\Omega} m(\bx) \dd\bx = \eta \,,
	\end{aligned}
\end{equation}
where the optimization problem is assumed to be wellposed and there exists a minimizer $\tilde{m}$. 

Next, let $\calF_{NN}$ is a neural operator approximation of $\calF$, and $J_{NN}(m)$ and $J^C_{NN}(m)$ are two approximations of $J(m)$ given by
\begin{equation*}
	J_{NN}(m) := \int_{\Gamma_{out}} g \, \calF_{NN}(m) \dS(\bx), \qquad J_{NN}^C(m) := \int_{\Gamma_{out}} g \, \calF^C(m, \calF_{NN}(m)) \dS(\bx)\,.
\end{equation*}
Let the minimizers of \eqref{eq:optProb} with the above two cost functions are denoted by $\tilde{m}_{NN}$ and $\tilde{m}^C_{NN}$, respectively. The main objective of this example is to compare the accuracy of $\tilde{m}_{NN}$ and $\tilde{m}^C_{NN}$ with $\tilde{m}$. 

\subsubsection{Data Generation, Neural Operators, and Numerical Method for the Optimization Problem}\label{sss:NNTopoOpt}

The domain $\Omega$ is triangulated using Gmsh\footnote{\url{https://gmsh.info/}} \cite{geuzaine2009gmsh} with 20614 triangular elements and 10301 vertices. The mesh is converted into a Fenics-friendly format using  Meshio\footnote{\url{https://github.com/nschloe/meshio}} \cite{schlomer_nico_2022_6346837}. The finite element spaces for the parameter and solution fields are based on the first-order Lagrange basis functions over the mesh of domain $\Omega$. If $\calM_h \subset \calM$ and $\calU_h \subset \calU$ are the finite element function spaces, then the variational problem in discrete setting reads:
\begin{equation}\label{eq:diffWeakProb2FE}
	\begin{aligned}
		\text{Given } m\in \calM_h, \text{ find }u \in \calU_h \text{ such that } \qquad\dualDot{v}{\calR(m, u)} = 0\,,\quad\forall v\in\calU_{h}\,.
	\end{aligned}
\end{equation}
Let, as before $\calF_h$, denote the finite-dimensional approximation of $\calF$. 

\paragraph{Data for Neural Operator Training}

\begin{figure}
	\begin{subfigure}[t]{.5\textwidth}
		\centering
		\includegraphics[width = 0.9\textwidth]{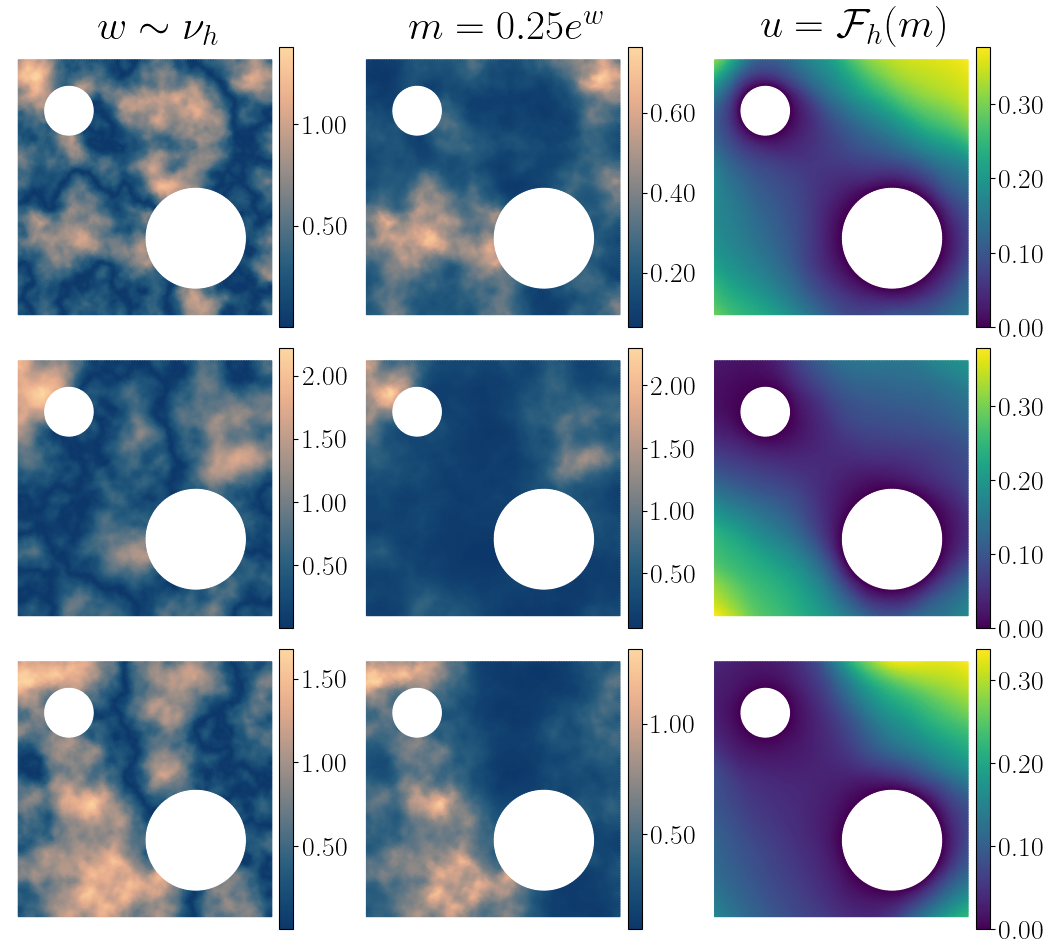}
		%\vspace{-10pt}
		\caption{}
	\end{subfigure}%
	\begin{subfigure}[t]{.5\textwidth}
		\centering
		\includegraphics[width = \textwidth]{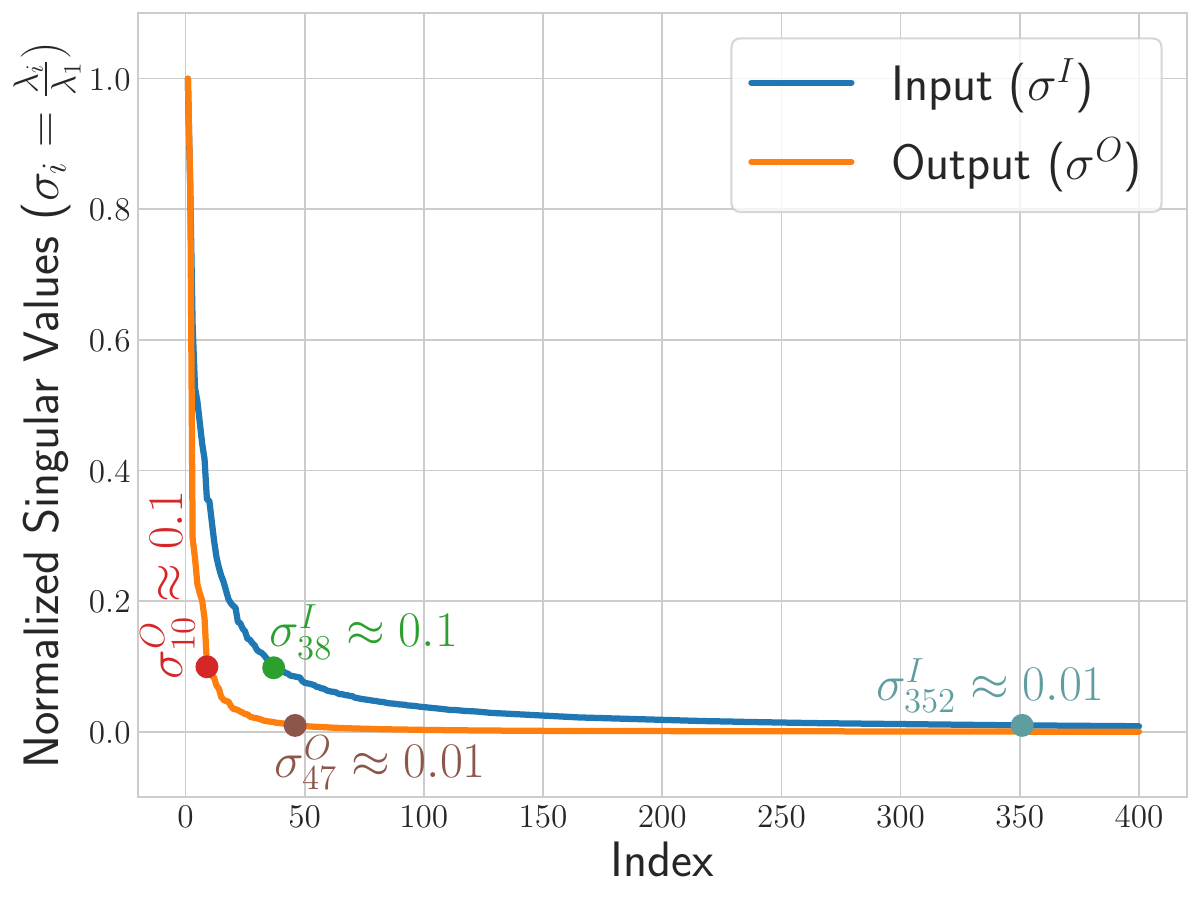}
		%\vspace{-10pt}
		\caption{}
	\end{subfigure}
	\vspace{-5pt}
	\caption{(a) Visualization of three representative samples, where $w$ and $m$ are from \eqref{eq:mDraw2} and $u$ is the solution of \eqref{eq:diffWeakProb2FE} given $m$. (b) Normalized singular values for input and output data of sample size $N = 4096$ for the second example. The indices associated with the singular values near $0.1$ and $0.01$ are also displayed. In (a), the first column displays samples $w$ generated from prior $\nu_h$, the second column shows the diffusivity $m = 0.25 e^{w^i(\bx)}$, and the last column is the finite element solution of the nonlinear reaction-diffusion equation with $m$ as diffusivity. We note that the pair $(m, u=\calF_h(m))$ constitutes a data point.}
	\label{fig:topoSamples}
\end{figure}

Let the probability distribution $\nu$ and its finite element approximation $\nu_h$ be the same as in \cref{sss:dataGen}. The training samples $\{(m^i, u^i)\}_{i=1}^N$ are generated as follows:
\begin{equation}\label{eq:mDraw2}
	\text{Draw } w^i(\bx) \sim \nu_h\, \qquad \text{ and } \qquad  m^i(\bx) = 0.25 e^{w^i(\bx)}\,, \quad u^i = \calF_h(m^i)\,.
\end{equation}
In contrast to the example in \cref{ss:prob1}, here, diffusivity $m = 0.25e^w$ is taken as the model parameter, where $w$ is some function. To generate training data, $w$ is sampled using the probability distribution $\nu_h$. In \cref{fig:topoSamples}(a), three representative data samples along with the singular values of input and output data from $N = 4096$ samples are depicted. 

\paragraph{Neural Operators}

\begin{figure}
	\centering
	\vspace{-5pt}
	\includegraphics[width = 0.6\textwidth]{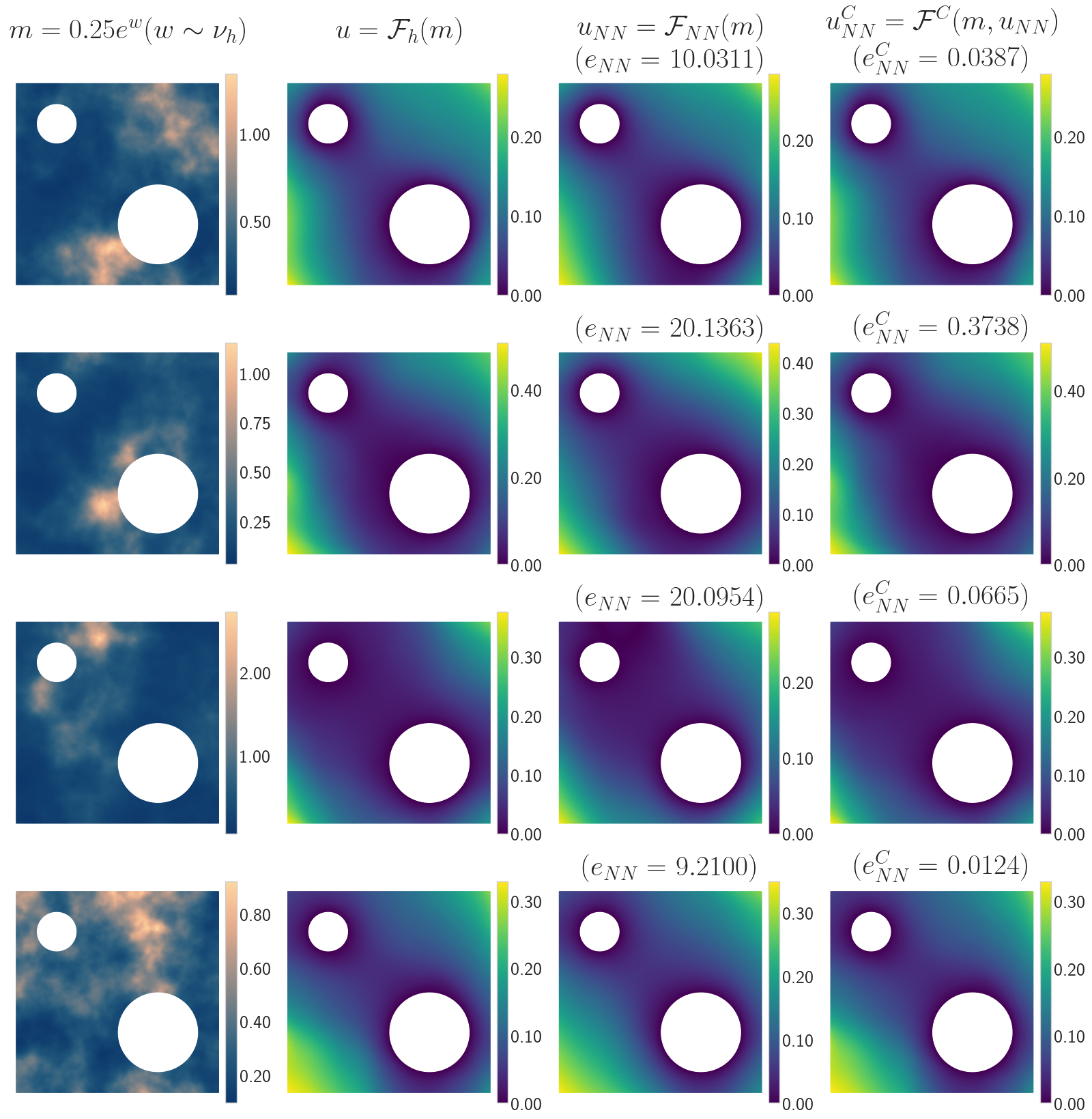}
	\vspace{-5pt}
	\caption{Comparing true solution, neural operator prediction, and the correction of neural operator prediction for networks $1, 6, 11, 16$ (see \cref{tab:topoCompareAcc}) trained with smaller samples of data.}
	\label{fig:compareAccVisOpt1}
\end{figure}

\begin{figure}
	\centering
	\vspace{-5pt}
	\includegraphics[width = 0.6\textwidth]{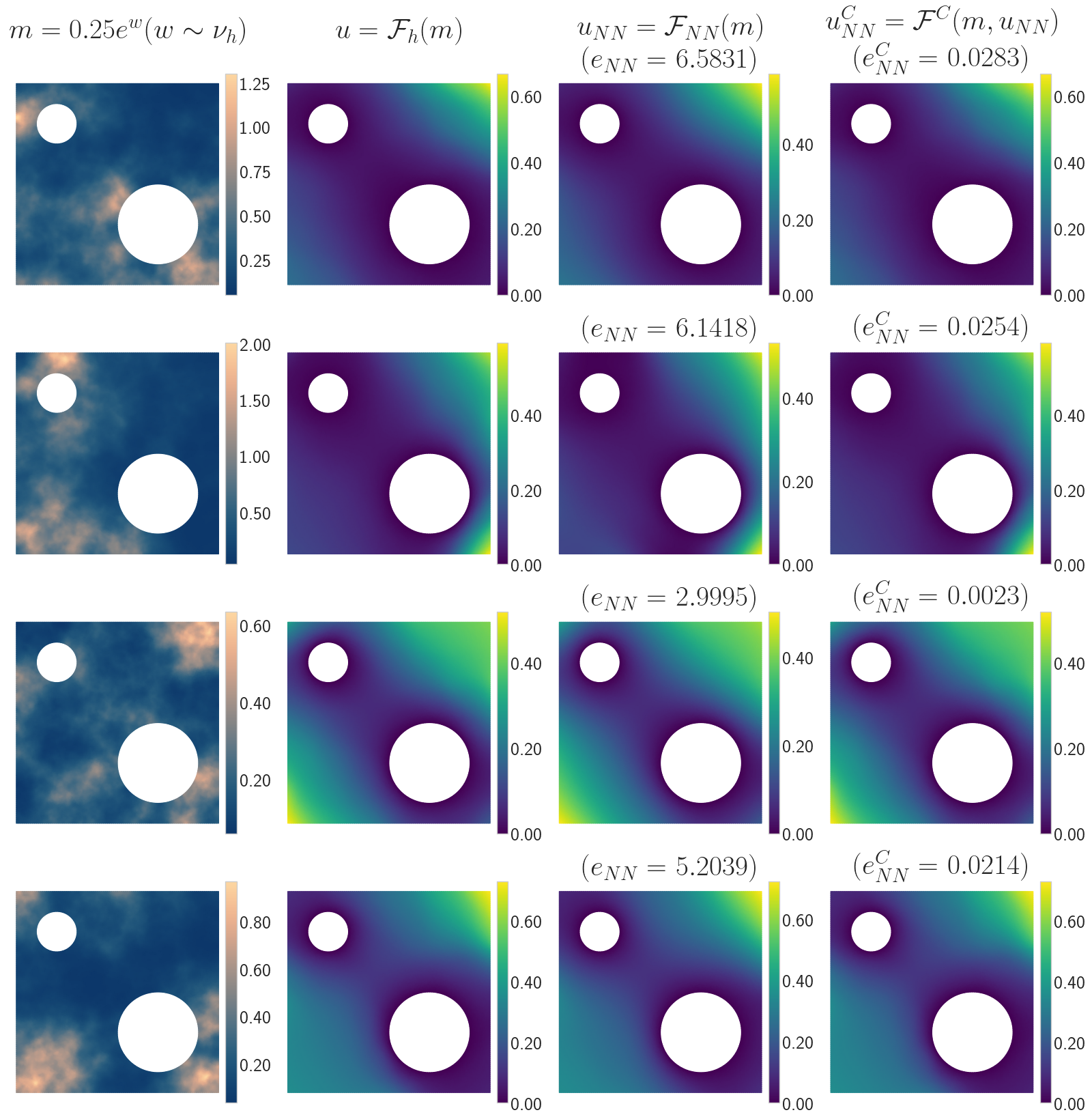}
	\vspace{-5pt}
	\caption{Comparing true solution, neural operator prediction, and the correction of neural operator prediction for networks $5, 10, 15, 20$ (see \cref{tab:topoCompareAcc}) trained with larger samples of data.}
	\label{fig:compareAccVisOpt2}
\end{figure}

Same as in the first example, a total of twenty neural operators with 
\begin{equation*}
	(r_m, r_u) \in \{(50, 25), (50, 50), (100, 25), (100, 50)\} \text{ and } N \in \{256, 512, 1024, 2048, 4096\}
\end{equation*}
 are trained. The plot of singular values in \cref{fig:topoSamples}(b) is similar to the first example as expected and shows that the singular values of output data decay faster relative to the input data. \cref{tab:topoCompareAcc} compares the accuracy of all twenty neural operators and corrections of these neural operators. Increasing data samples have a more dominant effect on the accuracy of neural operators in this example. It is however noted that increasing the reduced space dimension does not necessarily increase the accuracy. For the neural operators trained on the smallest and largest data samples, forward solutions and their approximations for a randomly drawn sample of $w$ with $m = 0.25e^w$ are compared in \cref{fig:compareAccVisOpt1} and \cref{fig:compareAccVisOpt2}.

{
\begin{table}
	\centering
	{\footnotesize
		%\rowcolors{3}{\tbcoli}{\tbcolii}
		\begin{tabular}{|c|c|c|c|c|c|>{\columncolor{\tbcoliii}}c|c|c|>{\columncolor{\tbcoliii}}c|}
			\hline
			&&&& \multicolumn{3}{c|}{$e_{NN}$} & \multicolumn{3}{c|}{$e^C_{NN}$} \\ \hline 
			Number & $r_m$ & $r_u$ & $N$ & $\mathrm{min}$ & $\mathrm{max}$ & $\mathrm{mean}$ & $\mathrm{min}$ & $\mathrm{max}$ & $\mathrm{mean}$ \\
			\hline
			\rowcolor{blue!10} 
			1 & 50 & 25 & 256 & 7.15830 & 21.25784 & \cellcolor{\tbcoliii} 12.47695 & 0.00153 & 0.84201 & \cellcolor{\tbcoliii} 0.13296 \\
			2 & 50 & 25 & 512 & 4.59103 & 15.99608 & 9.09779 & 0.00907 & 0.38084 & 0.06826 \\
			3 & 50 & 25 & 1024 & 4.56197 & 10.36562 & 6.26587 & 0.00459 & 0.16858 & 0.03795 \\
			4 & 50 & 25 & 2048 & 3.75125 & 7.91939 & 5.54681 & 0.00139 & 0.23467 & 0.04583 \\
			\rowcolor{blue!10} 
			5 & 50 & 25 & 4096 & 2.91057 & 8.09030 & \cellcolor{\tbcoliii} 4.94291 & 0.00178 & 0.12398 & \cellcolor{\tbcoliii} 0.02647 \\
			\rowcolor{green!10}
			6 & 50 & 50 & 256 & 9.15900 & 19.87726 & \cellcolor{\tbcoliii} 13.10120 & 0.02228 & 0.48271 & \cellcolor{\tbcoliii} 0.14268 \\
			7 & 50 & 50 & 512 & 5.60731 & 24.31812 & 11.78976 & 0.00593 & 0.77939 & 0.15389 \\
			8 & 50 & 50 & 1024 & 4.50428 & 16.09015 & 7.89234 & 0.00262 & 1.34911 & 0.10651 \\
			9 & 50 & 50 & 2048 & 3.58687 & 11.70707 & 6.41672 & 0.00263 & 0.29747 & 0.03527 \\
			\rowcolor{green!10}
			10 & 50 & 50 & 4096 & 3.11247 & 6.86457 & \cellcolor{\tbcoliii} 4.22486 & 0.00139 & 0.17169 & \cellcolor{\tbcoliii} 0.02497 \\
			\rowcolor{red!10}
			11 & 100 & 25 & 256 & 6.66220 & 22.61565 & \cellcolor{\tbcoliii} 11.11228 & 0.00673 & 1.41224 & \cellcolor{\tbcoliii} 0.20579 \\
			12 & 100 & 25 & 512 & 5.56698 & 21.51665 & 10.51670 & 0.00973 & 0.68513 & 0.09869 \\
			13 & 100 & 25 & 1024 & 5.25971 & 12.44923 & 7.97554 & 0.00135 & 0.59833 & 0.07555 \\
			14 & 100 & 25 & 2048 & 4.53138 & 11.74423 & 7.21405 & 0.00342 & 0.96087 & 0.11523 \\
			\rowcolor{red!10}
			15 & 100 & 25 & 4096 & 3.24184 & 6.81036 & \cellcolor{\tbcoliii} 4.38627 & 0.00031 & 0.18123 & \cellcolor{\tbcoliii} 0.02263 \\
			\rowcolor{orange!10}
			16 & 100 & 50 & 256 & 7.72910 & 18.16715 & \cellcolor{\tbcoliii} 11.80013 & 0.00926 & 0.33257 & \cellcolor{\tbcoliii} 0.09902 \\
			17 & 100 & 50 & 512 & 6.11802 & 12.84410 & 8.83076 & 0.00718 & 0.28607 & 0.05695 \\
			18 & 100 & 50 & 1024 & 4.99326 & 10.96458 & 7.80839 & 0.00267 & 0.25835 & 0.04400 \\
			19 & 100 & 50 & 2048 & 3.81907 & 15.93567 & 6.45772 & 0.00198 & 0.91764 & 0.06216 \\
			\rowcolor{orange!10}
			20 & 100 & 50 & 4096 & 3.41474 & 8.98157 & \cellcolor{\tbcoliii} 5.17495 & 0.00113 & 0.12705 & \cellcolor{\tbcoliii} 0.02403 \\
			\hline
		\end{tabular}
	}
	\caption{Comparing errors due to neural operator approximations of the forward model and corrections of neural operators for the second example. For more details on the table and colors in columns and rows, see \cref{tab:compareAcc}. Compared to the neural operators in the first example, the effect of increasing data samples on the average neural operator and corrector operator errors is more evident. However, increasing the input and output reduced dimensions have very little positive effect on the accuracy of neural operators.} \label{tab:topoCompareAcc}
\end{table}
}

\paragraph{Numerical Solution of the Optimization Problem}

The numerical method is based on the relaxation of the optimization problem \eqref{eq:optProb} in a finite-dimensional setting:
\begin{equation}\label{eq:optProbFE}
	\begin{aligned}
		&\qquad \min_{m \in \calM_h, \lambda\in \bbR} && \hat{J}(m, \lambda, u) := \int_{\Gamma_{out}} g\, u \dd \bx + \lambda \left( \int_{\Omega} m\dd \bx - \eta |\Omega| \right)\,, \\
		&\qquad \text{where } u = u(m)\in \calU_h \text{ satisfies }\qquad &&\dualDot{v}{\calR(m, u)} = 0\,, \qquad \forall v \in \calU_{h}, \\
		&\qquad \text{and }\qquad  &&0 < m_{\text{lw}} \leq m \leq 1\,.
	\end{aligned}
\end{equation}
By replacing $J_h(m)$ with 
\begin{equation*}
    J_{h, NN}(m) = \int_{\Gamma_{out}} g~\calF_{h, NN}(m) \dS(\bx) \qquad \text{ and } \qquad J^C_{h, NN}(m) = \int_{\Gamma_{out}} g~\calF^C(m, \calF_{h, NN}(m)) \dS(\bx)\,,
\end{equation*}
respectively, optimization problems with surrogates of the forward model are obtained. For the results in the next section, the numerical minimizers with cost functions $J_h$, $J_{h, NN}$, $J^C_{h, NN}$ are denoted by $\tilde{m}$, $\tilde{m}_{NN}$, and $\tilde{m}^C_{NN}$, respectively. 

Optimization problem \eqref{eq:optProbFE} is solved using a bi-level iterative scheme, wherein the outer iteration, $u$ and pair $(m, \lambda)$ are solved sequentially. For a given outer iteration step $k$ and variables $m_k, \lambda_k, u_k = u(m_k)$, first, the new updated values $(m_{k+1}, \lambda_{k+1})$ are computed using the inner iteration, and then using $m_{k+1}$, $u_{k+1} = u(m_{k+1})$ is computed. The outer iteration stops when $||m_k - m_{k+1}||_{L^2(\Omega)} \leq \gamma_{tol}$ or when $k = n_{max}$. The numerical method is detailed in \cref{s:topOptDetails}.

In the numerical experiments, the target average diffusivity is taken to be $\eta = 0.4$, and the initial guess for the parameter $m$ is a constant function $m(\bx) = 0.1$. Lagrange multiplier $\lambda$ is initialized as $\lambda = 1$. The tolerance in \cref{alg:optProbInner} is set to $m_{\text{tol}} = 0.005$.  Finally, $m_{\text{lw}} = 0.001$. 

\subsubsection{Optimization Results}\label{sss:optResults}

The optimization problem was solved numerically for the following three cases:
\begin{itemize}
	\item when $u = u(m)$ is obtained by solving the forward problem in which case the numerical minimizer is denoted by $\tilde{m}$;
	\item when neural operator prediction was used to approximate $u$ by $u_{NN} = \calF_{NN}(m)$ in which case the minimizers are identified by $\tilde{m}_{NN}$; and
	\item when $u$ was approximated by $u^C_{NN} = \calF^C(m, \calF_{NN}(m))$ using the corrector operator in which case the minimizer is denoted by $\tilde{m}^C_{NN}$. 
\end{itemize}

The percentage errors between numerical minimizers can be defined as:
\begin{equation}\label{eq:errsNormalizedM}
	\tilde{\varepsilon}_{NN} := \frac{||\tilde{m} - \tilde{m}_{NN}||_{l^2}}{||\tilde{m}||_{l^2}} \times 100\,, \qquad \qquad \tilde{\varepsilon}^C_{NN} := \frac{||\tilde{m} - \tilde{m}^C_{NN}||_{l^2}}{||\tilde{m}||_{l^2}} \times 100\,.
\end{equation}
Similarly, the errors in forward solutions when using the minimizers as the model parameter are defined as
\begin{equation}\label{eq:errsNormalizedU}
	\tilde{e}_{NN} := \frac{||u(\tilde{m}) - u_{NN}(\tilde{m}_{NN}) ||_{l^2}}{||u(\tilde{m})||_{l^2}} \times 100\,, \qquad \tilde{e}^C_{NN} := \frac{||u(\tilde{m}) - u^C_{NN}(\tilde{m}^C_{NN}) ||_{l^2}}{||u(\tilde{m})||_{l^2}} \times 100 \,.
\end{equation}

\begin{figure}
	\centering
	\vspace{-10pt}
	\includegraphics[width = 0.6\textwidth]{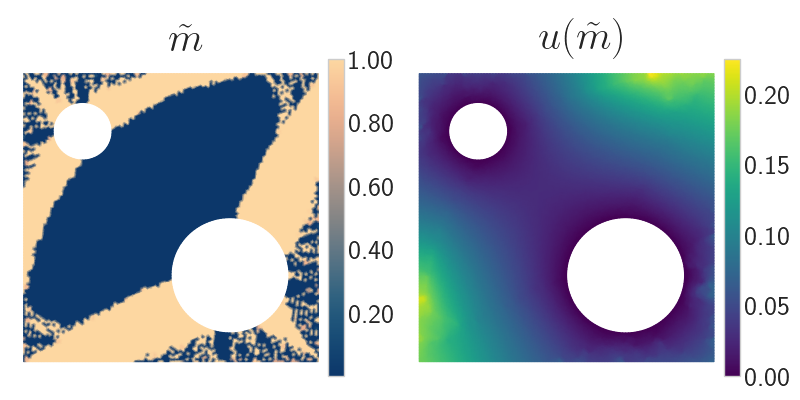}
	\vspace{-5pt}
	\caption{Plot of ``true" minimizer (up to numerical discretization error in forward and optimization problems) $\tilde{m}$ and corresponding forward solution $u(\tilde{m})$. }
	\label{fig:optResTrue}
\end{figure}

\begin{figure}
	\centering
	\vspace{-10pt}
	\includegraphics[width = 0.75\textwidth]{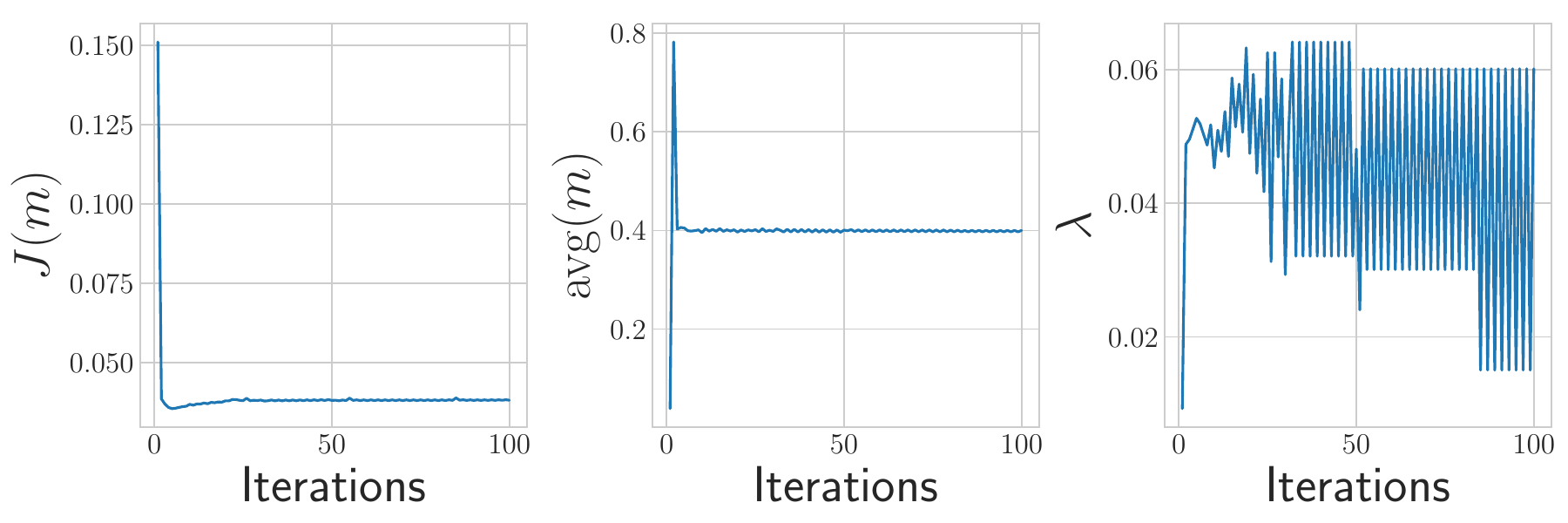}
	\caption{History of compliance function $J(m)$, the volume average of $m$ (note that the target volume average is $\eta = 0.4$), and the Lagrange multiplier during iterations.}
	\vspace{-5pt}
	\label{fig:optResTrueConvg}
\end{figure}

\begin{figure}
	\centering
	\vspace{-10pt}
	\includegraphics[width = 0.8\textwidth]{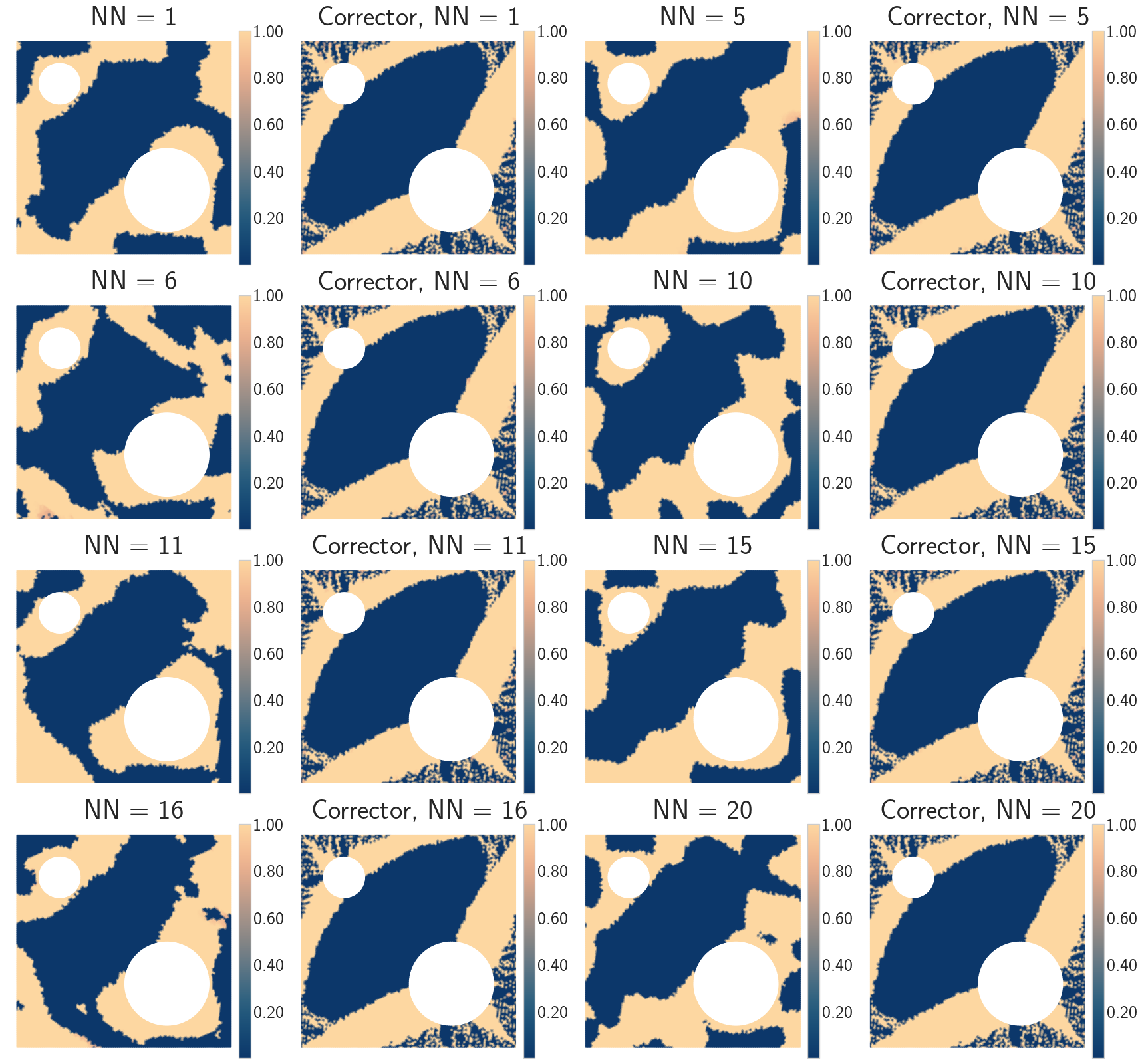}
	\vspace{-5pt}
	\caption{Comparing the minimizers for neural networks trained with a smaller dataset (left two columns) and the larger dataset (right two columns). The odd columns (1 and 3) correspond to the minimizers obtained by employing neural operator surrogates in the cost function. In contrast, the even columns (2 and 4) are those where neural operator predictions are corrected using the corrector operator $\calF^C$. For the properties of neural operators, see to \cref{tab:topoCompareAcc}.}
	\label{fig:optResVis}
\end{figure}

\begin{figure}
	\centering
	%\vspace{-10pt}
	\includegraphics[width = 0.75\textwidth]{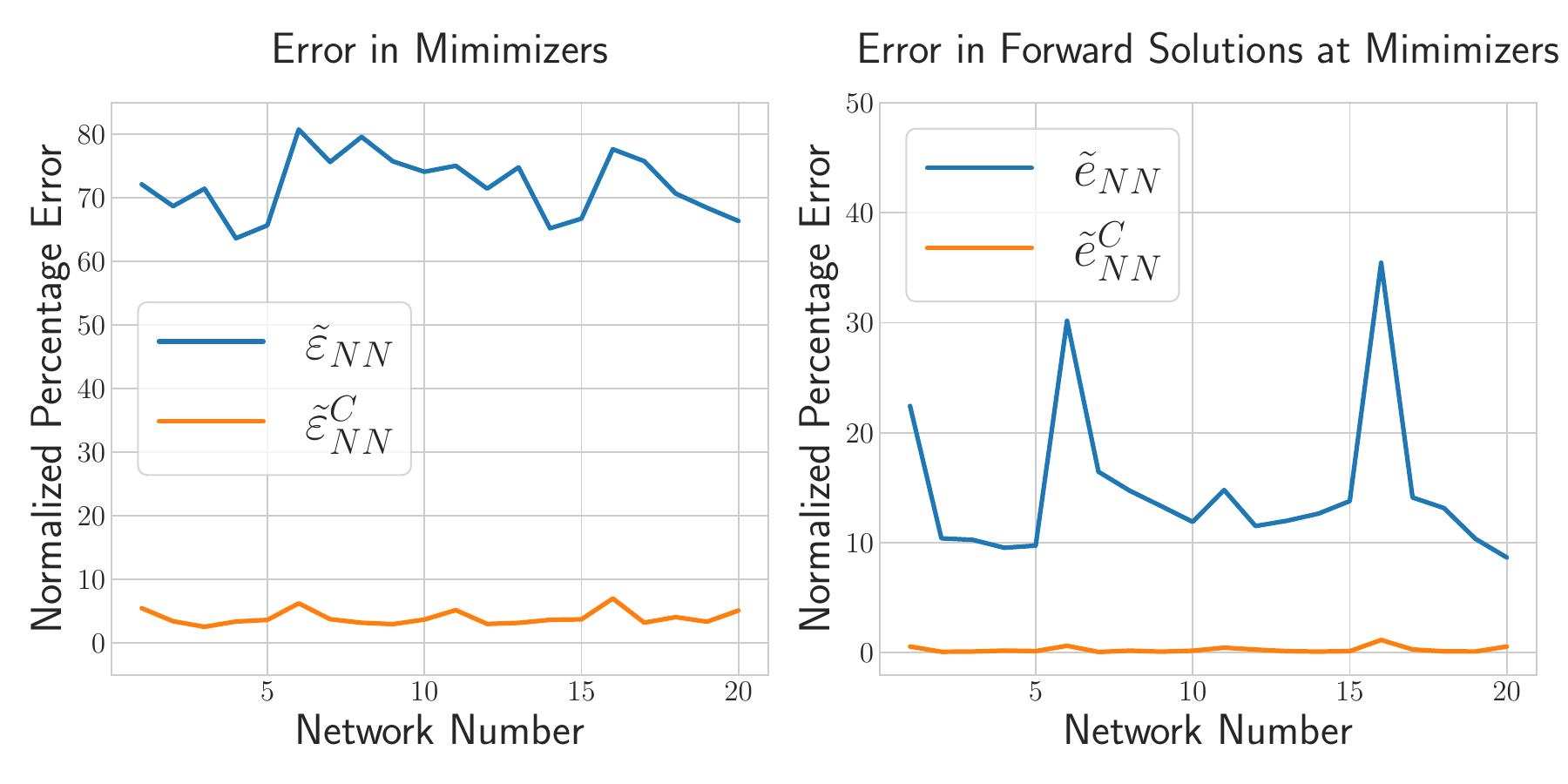}
	%\vspace{-10pt}
	\caption{Normalized minimizer errors $\tilde{\varepsilon}_{NN}$ and $\tilde{\varepsilon}_{NN}^C$ due to surrogate approximations of the forward problem (see \eqref{eq:errsNormalizedM}) and the error in forward solutions $\tilde{e}_{NN}$ and $\tilde{e}^C_{NN}$ defined in \eqref{eq:errsNormalizedU}.}
	\label{fig:optResl2Norm}
\end{figure}

The numerical minimizer $\tilde{m}$ and the forward solution at $\tilde{m}$ are shown in \cref{fig:optResTrue}. The history of the cost function, the volumetric average of $m$, and the Lagrange multiplier are plotted in \cref{fig:optResTrueConvg}. In \cref{fig:optResVis}, the minimizers for the neural operators trained with smaller and larger datasets are compared. The figure also shows the minimizers when neural operators are corrected using the corrector operator $\calF^C$. Finally, for all the twenty neural operators and the corrector of those neural operators, the percentage errors $\tilde{\varepsilon}_{NN}$, $\tilde{\varepsilon}^C_{NN}$, $\tilde{e}_{NN}$, and $\tilde{e}^C_{NN}$ are plotted in \cref{fig:optResl2Norm}. The errors are tabulated in \cref{tab:optResl2Norm} to allow easier comparison of the accuracy of neural operators with and without corrections. From the error results in \cref{tab:optResl2Norm}, it is clear that the neural operators consistently lead to minimizers with high error (as high as 80 percent). The corrector on the other hand provides an approximation of minimizers with errors below seven percent. 

{
	\begin{table}
		\centering
		{\footnotesize
			%\rowcolors{3}{\tbcoli}{\tbcolii}
			\begin{tabular}{|c|c|c|c|c|c|c|c|c|c|}
				\hline
				Number & $\tilde{\varepsilon}_{NN}$ & $\tilde{\varepsilon}^C_{NN}$ & $\tilde{e}_{NN}$ & $\tilde{e}^C_{NN}$ & Number & $\tilde{\varepsilon}_{NN}$ & $\tilde{\varepsilon}^C_{NN}$ & $\tilde{e}_{NN}$ & $\tilde{e}^C_{NN}$ \\
				\hline
				\rowcolor{blue!10} 
				1 & 72.14426 & 5.46892 & 22.43047 & 0.56382 & 11 & 75.05329 & 5.16826 & 14.79259 & 0.45834 \\
				2 & 68.70190 & 3.41497 & 10.39295 & 0.08882 & 12 & 71.45512 & 2.99042 & 11.51458 & 0.28541 \\
				3 & 71.45517 & 2.53612 & 10.25896 & 0.10772 & 13 & 74.81585 & 3.15803 & 11.99392 & 0.14022 \\
				4 & 63.65149 & 3.36944 & 9.54038 & 0.19412 & 14 & 65.21906 & 3.63778 & 12.64035 & 0.10310 \\
				\rowcolor{green!10}
				5 & 65.68216 & 3.62392 & 9.73351 & 0.14817 & 15 & 66.73859 & 3.70730 & 13.78046 & 0.14303 \\
				\rowcolor{blue!10}
				6 & 80.76804 & 6.22848 & 30.17278 & 0.63044 & 16 & 77.65764 & 6.99276 & 35.45096 & 1.16036 \\
				7 & 75.65322 & 3.72890 & 16.44920 & 0.07026 & 17 & 75.78578 & 3.19518 & 14.10578 & 0.29605 \\
				8 & 79.60037 & 3.17328 & 14.71594 & 0.18724 & 18 & 70.69823 & 4.06247 & 13.14865 & 0.12392 \\
				9 & 75.74455 & 2.96135 & 13.32289 & 0.09598 & 19 & 68.44559 & 3.34363 & 10.35292 & 0.11209 \\
				\rowcolor{green!10}
				10 & 74.11825 & 3.67819 & 11.89410 & 0.18198 & 20 & 66.36491 & 5.09965 & 8.65368 & 0.56261 \\
				\hline
			\end{tabular}
		}
		\caption{Comparison of the percentage errors in minimizers and forward solutions at the minimizers for neural operators and corrected neural operators. Networks trained with small datasets are highlighted in blue, while the ones trained with relatively large datasets are highlighted in green.} \label{tab:optResl2Norm}
	\end{table}
}

%%%% --------------------------- %%%%
\section{Conclusion}\label{s:conc}

The work considers a powerful approach for enhancing the accuracy and reliability of neural operators, especially when the neural operator accuracy is impacted by the unavailability of appropriate training distributions and sparse data. The approach is based on the corrector operator, which requires solving the linear variational problem given the input parameter $m$ and the prediction furnished by the neural operator. For the two examples considered in this work, the increase in accuracy obtained via the corrector operator is not possible to attain by simply tuning the hyperparameters of the neural operators and increasing the training data samples. For the topology optimization problem, the work highlights the limitations of neural operators in choosing the appropriate training samples and neural operators leading to results with large errors. The corrector operator for the topology optimization problem increased the accuracy of optimizers significantly. In summary, the approach seems to do a great job of increasing the accuracy and reliability of neural operators. 

In the present work, the correction is computed external to neural operators. In contrast, future work will explore the possibility of integrating the correction step into the neural operators. Further, goal-oriented error estimates can be used to enhance the accuracy of neural operators with respect to specific quantities of interest. Finally, downstream applications of parameter estimation and design optimization of complex materials for mechanical loading and actuation are of interest where the neural operators will be used as surrogates of highly nonlinear parametric multiphysics models of mechanical deformation. 

%%%% --------------------------- %%%%
\section*{Acknowledgement}
The author is thankful to the Late Professor John Tinsley Oden for editing and making corrections in several versions of the draft before it could be ready for publication. Professor Oden also helped write part of the second section of the article for which the author is grateful. This work was supported by the U.S. Department of Energy, Office of Science, USA, Office of Advanced Scientific Computing Research, Mathematical Multifaceted Integrated Capability Centers (MMICCS), under Award Number DE-SC0019303.

%-------------
% References
%-------------
\addcontentsline{toc}{section}{References}
%\bibliographystyle{model1-num-names}
%\biboptions{sort,numbers,comma,compress}                 
%\bibliography{main.bib}

%%%% --------------------------- %%%%
\appendix

\section{Corrector Operator Error Analysis}\label{s:corrErrAnalysis}

In this section, \cref{thm:corrOpConvergence} is proved. Let $m \in \calM$, and $\tilde{u} \in \calU$ be arbitrary approximation of $u = \calF(m)$, where $\calF$ is the solution operator. From the definition of corrector operator $\calF^C$ (see \eqref{eq:corrOp}), it can be shown that
\begin{equation}\label{eq:fcAnalysis1}
\begin{split}
    &\calF^C(m) = u^C = \tilde{u} - \delta_u\calR(m, \tilde{u})^{-1} \calR(m, \tilde{u}) \\
    \Rightarrow \qquad & \underbrace{u - u^C}_{= \, e^C} = \underbrace{u - \tilde{u}}_{= \, \tilde{e}} + \delta_u \calR(m, \tilde{u})^{-1}\left(\calR(m, \tilde{u}) - \underbrace{\calR(m, u)}_{= \, 0}\right) \\
    \Rightarrow \qquad & e^C = \tilde{e} - \delta_u \calR(m, \tilde{u})^{-1} \left(\calR(m, u) - \calR(m, \tilde{u})\right) \,.
\end{split}
\end{equation}
Using the identity $\tilde{e} = \delta_{u} \calR(m, \tilde{u})^{-1} \delta_{u}\calR(m, \tilde{u})(\tilde{e})$ and the Taylor series expansion
\begin{equation*}
    \calR(m, u) - \calR(m, \tilde{u}) = \delta_u \calR(m, \tilde{u})(\tilde{e}) + \int_0^1 (1 - s) \delta_u^2 \calR(m, \tilde{u} + s\tilde{e})(\tilde{e}, \tilde{e}) \dd s
\end{equation*}
in \eqref{eq:fcAnalysis1}, it holds
\begin{equation}\label{eq:fcAnalysis2}
    \begin{split}
        e^C &=  \delta_{u} \calR(m, \tilde{u})^{-1} \delta_{u}\calR(m, \tilde{u})(\tilde{e}) - \delta_u \calR(m, \tilde{u})^{-1} \left(\calR(m, u) - \calR(m, \tilde{u})\right) \\
        &= \delta_{u} \calR(m, \tilde{u})^{-1} \left[  \delta_{u}\calR(m, \tilde{u})(\tilde{e}) - \left(\calR(m, u) - \calR(m, \tilde{u})\right) \right] \\
        &= \delta_{u} \calR(m, \tilde{u})^{-1} \left[ - \int_0^1 (1-s) \delta_u^2 \calR(m, \tilde{u} + s\tilde{e})(\tilde{e}, \tilde{e}) \dd s \right] \\
        &= \int_0^1 (s-1) \delta_{u} \calR(m, \tilde{u})^{-1} \delta_u^2 \calR(m, \tilde{u} + s\tilde{e})(\tilde{e}, \tilde{e}) \dd s  \,.
    \end{split}
\end{equation}
Next, define a bilinear operator $A(s): \calU \times \calU \to \calU$ for $s\in [0,1]$ as follows
\begin{equation*}
    A(s)(p, q) := \delta_{u} \calR(m, \tilde{u})^{-1} \delta_u^2 \calR(m, \tilde{u} + s\tilde{e})(p, q)
\end{equation*}
and note that
\begin{equation*}
    ||A(s)(p, q)||_{\calU} \leq ||A(s)||_{\calL(\calU \times \calU; \calU)} \, ||p||_{\calU} \, ||q||_{\calU} \leq \left[\sup_{s\in [0, 1]} ||A(s)||_{\calL(\calU \times \calU; \calU)} \right] \, ||p||_{\calU} \, ||q||_{\calU}\,.
\end{equation*}
Taking the norm of both sides in \eqref{eq:fcAnalysis2} and using the above estimate, the following can be shown
\begin{equation}
    ||e^C||_{\calU} \leq \int_0^1 (1-s) ||A(s)(\tilde{e}, \tilde{e})||_{\calU} \leq  \frac{1}{2} \left[\sup_{s\in [0, 1]} \vert\vert \delta_{u} \calR(m, \tilde{u})^{-1} \delta_u^2 \calR(m, \tilde{u} + s\tilde{e}) \vert\vert_{\calL(\calU \times \calU; \calU)} \right] \, ||\tilde{e}||_{\calU}^2\,.
\end{equation}
\qed

\section{Bounds on Derivatives of a Residual for Nonlinear Diffusion Example}\label{s:bdResDerNonDiff}

In this section, \cref{thm:bdResDerNonDiff} is established. In what follows, the preliminary results needed in the proof are first collected, and then \cref{thm:bdResDerNonDiff} is proved.

\subsection{Preliminary Results}
Let $X$ and $Y$ be two Banach spaces. $X$ is said to be embedded continuously in $Y$, written as $X \hookrightarrow Y$, if 
\begin{itemize}
    \item $X \subseteq Y$;
    \item the canonical injection $i: X \to Y$ is a continuous (linear) operator, i.e., there exists a constant $C>0$ such that
    \begin{equation*}
        ||i(u)||_Y \leq  C ||u||_X, \qquad \forall u \in X\,.
    \end{equation*}
    Here, $C$ is independent of $u$.
\end{itemize}
Embedding of $X$ into $Y$ is compact if, in addition to the above conditions, the canonical injection operator $i$ is a compact operator. Relevant Sobolev embedding results from Theorems 2.6.1 and 2.6.2 in \cite{badiale2010semilinear} are collected in the following theorem:
\begin{theorem}
Let $\Omega \subset \R^d$ be the open, bounded, and smooth domain where $d = 2, 3$. It holds that
\begin{equation*}
    H^1(\Omega) \hookrightarrow L^q(\Omega)\,,
\end{equation*}
for every $q\in [1, \infty)$ when $d=2$ and for every $q \in [1, 6]$ when $d=3$. Further, the embedding is compact for every $q\in [1, \infty)$ when $d=2$ and for every $q \in [1, 6)$ when $d=3$.

For $d=3$, if $\Omega$ is open and bounded subset of $\R^d$ or $\Omega = \R^d$, then
\begin{equation*}
    H^1_0(\Omega) \hookrightarrow L^q(\Omega)
\end{equation*}
for every $d \in [1, 6]$. Further, the embedding is compact for every $q\in [1, 6)$. 
\end{theorem}

Consider the case when $\Omega$ is open, bounded, and smooth with $d=2$ or $d=3$. Let $\calU := H^1_0(\Omega)$. Then for every $u\in \calU$, using the above Sobolev embedding theorem, it holds
\begin{equation}\label{eq:sobEmb}
    ||u||_{L^4(\Omega)} \leq C_S ||u||_{\calU}\,,
\end{equation}
where $C_S$ independent of $u$ is the Sobolev embedding constant. 

Next, the Poincar\'e inequality theorem following [Theorem 2.6.3, \cite{badiale2010semilinear}] is stated:
\begin{theorem}
    Let $\Omega \subset \R^d$ be open and bounded. Then there exists a constant $C_P > 0$, depending only on $\Omega$, such that
    \begin{equation}\label{eq:poincare}
        ||u||_{L^2(\Omega)} \leq C_P ||\nabla u||_{L^2(\Omega)}, \qquad \forall u \in H^1_0(\Omega)\,.
    \end{equation}
\end{theorem}

\subsection{Proof of \cref{thm:bdResDerNonDiff}}

The proof is divided into three steps as follows.

\textbf{Step 1.} First, the upper bound on $||\delta_u \calR(m, \tilde{u})(v)||_{\calU^\ast}$ is established. Using the definition of operator norm gives
\begin{equation}\label{eq:bdEst1}
	||\delta_u \calR(m, \tilde{u}) ||_{\calL(\calU; \calU^\ast)} = \sup_{||v||_{\calU} = 1} ||\delta_u \calR(m, \tilde{u})(v)||_{\calU^\ast}\,.
\end{equation}
Further,
\begin{equation}\label{eq:bdEst2}
	||\delta_u \calR(m, \tilde{u})(v)||_{\calU^\ast} = \sup_{||w||_{\calU} = 1} \vert \langle w, \delta_u \calR(m, \tilde{u})(v) \rangle_{\calU} \vert \,.
\end{equation}
Noting the expression of $\delta_{u} \calR(m, u)$ in \eqref{eq:diffExampleResDer} and setup of the problem described in \cref{ss:nonDiff}, it holds
\begin{equation*}
    \begin{split}
        \vert \langle w, \delta_u \calR(m, \tilde{u})(v) \rangle_{\calU} \vert &\leq \int_\Omega \left\{\left( \kappa_0 \sup_{x\in \Omega} |m(x)| \right) |\nabla v| \,|\nabla w| + 3 \alpha \tilde{u}^2 v w \right\} \dd x \\
        &\leq \left[ \kappa_0 ||m||_{L^\infty(\Omega)} ||\nabla v||_{L^2(\Omega)} \, ||\nabla w||_{L^2(\Omega)} + 3 \alpha \left( \int_\Omega \tilde{u}^4 \dd x\right)^{1/2} \, \left( \int_\Omega v^2 w^2 \dd x\right)^{1/2} \right] \\
        &\leq \kappa_0 ||m||_{L^\infty(\Omega)} ||v||_{\calU} \, ||w||_{\calU} + 3 \alpha ||\tilde{u}||_{L^4(\Omega)}^2 \, ||v||_{L^4(\Omega)}\, ||w||_{L^4(\Omega)}\\
        &\leq \kappa_0 ||m||_{L^\infty(\Omega)} ||v||_{\calU} \, ||w||_{\calU} + 3 \alpha C_S^4 ||\tilde{u}||_{\calU}^2 \, ||v||_{\calU}\, ||w||_{\calU}\,,
    \end{split}
\end{equation*}
where H\"older inequality is used in the second and third equations, and the Sobolev embedding property (see \eqref{eq:sobEmb}) in the last equation. Using the above estimate in \eqref{eq:bdEst2} and combining the result with \eqref{eq:bdEst1}, it can be shown that 
\begin{equation}\label{eq:bdEst2.1}
	||\delta_u \calR(m, \tilde{u})(v)||_{\calU^\ast} \leq \underbrace{\left(\kappa_0 ||m||_{L^\infty(\Omega)} + 3 \alpha C_S^4 ||\tilde{u}||_{\calU}^2 \right)}_{\bar{C}_{\delta R}} \, ||v||_{\calU} \,.
\end{equation}

For the lower bound on $||\delta_u \calR(m, \tilde{u})(v)||_{\calU^\ast}$, proceeding as follows
\begin{equation*}
    \begin{split}
        ||\delta_u \calR(m, \tilde{u})(v)||_{\calU^\ast} &= \sup_{||w||_{\calU} = 1} \vert \langle w, \delta_u \calR(m, \tilde{u})(v) \rangle_{\calU} \vert \\
        &\geq \frac{\vert \langle v, \delta_u \calR(m, \tilde{u})(v) \rangle_{\calU} \vert}{||v||_{\calU}} \\
        &= \frac{1}{||v||_{\calU}}\left\vert \int_{\Omega} \left\{ \kappa_0 m(x) |\nabla v|^2 + 3 \alpha \tilde{u}^2 v^2 \right\} \dd x \right\vert \\
        &\geq \frac{\kappa_0 m_{\text{lw}}}{||v||_{\calU}} ||\nabla v||_{L^2(\Omega)}^2 = \frac{\kappa_0 m_{\text{lw}}}{2||v||_{\calU}} \left[ ||\nabla v||_{L^2(\Omega)}^2 + ||\nabla v||_{L^2(\Omega)}^2 \right] \\
        &\geq \frac{\kappa_0 m_{\text{lw}}}{2||v||_{\calU}} \left[ ||\nabla v||_{L^2(\Omega)}^2 + C_P^{-2} ||v||_{L^2(\Omega)}^2 \right] \\
        &\geq \underbrace{\frac{\kappa_0 m_{\text{lw}}}{2} \min\{1, C_P^{-2}\}}_{=:\hat{C}_{\delta R}} \, \frac{1}{||v||_{\calU}} \, ||v||_{\calU}^2 = \hat{C}_{\delta R}\, ||v||_{\calU}\,,
    \end{split}
\end{equation*}
where $\sup_{||w||_{\calU} = 1}  \vert \langle w, \delta_u \calR(m, \tilde{u})(v) \rangle_{\calU} \vert \geq \frac{\vert \langle v, \delta_u \calR(m, \tilde{u})(v) \rangle_{\calU} \vert}{||v||_{\calU}}$ is used in the second equation. In the fourth equation, the positive term is dropped, and the property that $m \geq m_{\text{lw}}$ is used, since $m \in \calM = \{g\in L^2(\Omega)\cap L^\infty(\Omega): g \geq m_{\text{lw}}\}$. In the fifth step, the Poincar\'e inequality from \eqref{eq:poincare} is applied. Combining the above with \eqref{eq:bdEst2.1}, the proof of (i) of \cref{thm:bdResDerNonDiff} is complete. 

\textbf{Step 2.} In this step, the upper bound on the norm of $\delta_u \calR(m, \tilde{u})^{-1}$ is established. Let $A:= \delta_u \calR(m, \tilde{u})$ and $r = A(v) \in \calU^\ast$, then
\begin{equation*}
    ||A^{-1}(r)||_{\calU} = ||A^{-1}(A(v))||_{\calU} = ||v||_{\calU} \,. 
\end{equation*}
However, since $r = A(v)$, from the previous calculations, it holds that $\hat{C}_{\delta R} ||v||_{\calU} \leq ||r||_{\calU^\ast}$, i.e., $||v||_{\calU} \leq \frac{1}{\hat{C}_{\delta R}} ||r||_{\calU^\ast}$. Combining this with the equation above, the following can be shown
\begin{equation*}
    ||A^{-1}(r)||_{\calU} = ||v||_{\calU} \leq {\hat{C}_{\delta R}} ||r||_{\calU^\ast}\,.
\end{equation*}
The above establishes (ii) of \cref{thm:bdResDerNonDiff}.

\textbf{Step 3.} To prove the last result, i.e., (iii) of \cref{thm:bdResDerNonDiff}, the definition of an operator norm is expanded as follows
\begin{equation*}
    ||\delta_u^2 \calR(m, \tilde{u}) ||_{\calL(\calU\times \calU; \calU^\ast)} = \sup_{||v||_{\calU} = 1, ||w||_{\calU} = 1} ||\delta_u^2 \calR(m, \tilde{u})(v,w)||_{\calU^\ast}\,.
\end{equation*}
Focusing on the $||\delta_u^2 \calR(m, \tilde{u})(v,w)||_{\calU^\ast}$, it holds that
\begin{equation}\label{eq:bdEst3}
\begin{split}
    ||\delta_u^2 \calR(m, \tilde{u})(v,w)||_{\calU^\ast} &= \sup_{||q||_{\calU} = 1} \vert \langle q, \delta_u^2 \calR(m, \tilde{u})(v,w) \rangle_{\calU} \vert \\
    &= \sup_{||q||_{\calU} = 1} \left\vert \int_\Omega 6 \alpha \tilde{u} v w q \dd x \right\vert \leq 6 \alpha \sup_{||q||_{\calU} = 1} \, ||\tilde{u}||_{L^4(\Omega)}\, ||v||_{L^4(\Omega)}\, ||w||_{L^4(\Omega)}\, ||q||_{L^4(\Omega)} \\
    &\leq 6\alpha C_S^4 \sup_{||q||_{\calU} = 1} \,||\tilde{u}||_{\calU}\, ||v||_{\calU}\, ||w||_{\calU}\, ||q||_{\calU} \\
    &= \underbrace{\left(6 \alpha C_S^4 ||\tilde{u}||_{\calU}\right)}_{=: \bar{C}_{\delta^2 R}}\, ||v||_{\calU}\, ||w||_{\calU}\,,
\end{split}
\end{equation}
where H\"older inequality is applied twice, and the Sobolev embedding property is used. This completes the proof of \cref{thm:bdResDerNonDiff}.

\section{Numerical Method for the Topology Optimization Problem}\label{s:topOptDetails}

A bi-level iteration scheme is developed for solving the optimization problem \eqref{eq:optProbFE}. The method employed here is based on the so-called SIMP (Solid Isotropic Material with Penalization) method and closely follows the topology optimization example\footnote{\url{https://comet-fenics.readthedocs.io/en/latest/demo/topology_optimization/simp_topology_optimization.html}} in \cite{bleyer2018numericaltours}.

Since $u = u(m) \in \calU$ solves $\dualDot{v}{\calR(m, u)} = 0$ for all $v\in \calU$, taking $v = u$, and using the definition of $\calR$ from \eqref{eq:diffWeakProb2}, it can be easily shown that
\begin{equation*}
	\int_{\Gamma_{out}} g \, u(m) \dS(\bx) = \int_{\Omega} m\nabla u\cdot \nabla u \dd \bx + \int_{\Omega} u^4 \dd \bx\,.
\end{equation*}
Let $\bq = \bq(m) = m \nabla u(m)$ denote the flux and $e = e(m) = \bq(m) \cdot \bq(m)$, then from the above, the compliance can be expressed as
\begin{equation*}
	\int_{\Gamma_{out}} g \, u(m) \dS(\bx) = \int_{\Omega} \left[ \frac{e(m)}{m} + u(m)^4 \right] \dd \bx\,.
\end{equation*}
Using the relation above, the optimization problem \eqref{eq:optProbFE} can be re-written as follows
\begin{equation}\label{eq:optProbFE2}
	\begin{aligned}
		& \min_{m \in \calM_h, \lambda\in \bbR} && \hat{J}(m, \lambda, u) = \int_{\Omega} \left[ \frac{e(m)}{m} + u(m)^4 \right] \dd \bx + \lambda \left( \int_{\Omega} m\dd \bx - \eta |\Omega| \right)\,, \\
		& \text{where } u = u(m)\in \calU_h \text{ satisfies }\qquad &&\dualDot{v}{\calR(m, u)} = 0\,, \qquad \forall v \in \calU_{h}, \\
		& \text{and }\qquad  &&0 < m_{\text{lw}} \leq m \leq 1\,.
	\end{aligned}
\end{equation}

The problem above is coupled in variables $m$, $\lambda$, and $u$. To simplify the computation, an iterative scheme is sought in which the variables can be uncoupled and sequentially computed. One possible approach is to consider an iteration (referring to this as outer iteration) where given solutions from the previous iteration, i.e., $m_k, \lambda_k, u_k = u(m_k)$, first, the pair $(m_{k+1}, \lambda_{k+1})$ is updated by solving the approximate minimization problem with fixed $u = u_k$, i.e.,
\begin{equation}\label{eq:optProbFE3}
	\begin{aligned}
		&\qquad (m_{k+1}, \lambda_{k+1}) =  \argmin_{m \in \calM_h, \lambda\in \bbR} && \hat{J}(m, \lambda, u_k) = \int_{\Omega} \left[ \frac{e(m_k)}{m} + u(m_k)^4 \right] \dd \bx + \lambda \left( \int_{\Omega} m\dd \bx - \eta |\Omega| \right)\,, \\
		&\qquad \text{and }\qquad  &&0 < m_{\text{lw}} \leq m \leq 1\,.
	\end{aligned}
\end{equation}
Next, given $m_{k+1}$, $\calR(m_{k+1}, u_{k+1}) = 0$ is solved for $u_{k+1}$. The iteration over $k$ is continued until $||m_{k+1} - m_k||_{L^2(\Omega)} \leq \gamma_{tol}$ or until $k$ reaches the maximum number of iterations. The algorithm for this outer iteration is presented in \cref{alg:optProbOuter}. 

\begin{algorithm}[H]
	\caption{Outer iteration for solving \eqref{eq:optProbFE} or equivalently \eqref{eq:optProbFE2}.
	}
	\label{alg:optProbOuter}
		\SetAlgoLined
		{\bf Input:}\\
		(a) set up mesh, variational problem, and optimization parameters \; 
		(b) take initial guess $m_0$, $\lambda_0 = 1$\;
		{\bf Outer iteration:}\\
		$k = 0$\;
		\While{$k<n^{out}_{\text{max iter}}$}{
			(1) given $m_k$, $\lambda_k$, and $u_k$, solve for $(m_{k+1}, \lambda_{k+1}) \in \calM_h \times \bbR$ following \cref{alg:optProbInner} \;
			(2) solve $\dualDot{v}{\calR(m_{k+1}, u_{k+1})} = 0$, $\forall v\in \calU_h$, for $u_{k+1} \in \calU_h$ \;
			(3) update flux $\bq_{k+1} = m_{k+1} \nabla u_{k+1}$ and $e_{k+1} = \bq_{k+1} \cdot \bq_{k+1}$ \;
			(4) $m_{k} = m_{k+1}$, $\lambda_{k} = \lambda_{k+1}$, $u_k = u_{k+1}$, and $k \leftarrow k+1$ \;
			\lIf{$||m_{k+1} - m_{k}||_{L^2(\Omega)} < \gamma_{tol}$}
			{
				{\bf break}
			}
		}
		{\bf Return: } $\tilde{m} = m_{k}$.
\end{algorithm}

Focusing now on \eqref{eq:optProbFE3}, the problem is solved iteratively (inner iteration) where in each iteration $i$, first $m^{(i)}$ is updated to compute $m^{(i+1)}$ while keeping $\lambda^{(i)}$ fixed and then $\lambda^{(i)}$ is updated into $\lambda^{(i+1)}$. During the inner iteration inside the outer iteration step $k$, $u$ is fixed to $u_k = u(m_k)$ throughout; \cref{alg:optProbOuter}. The key question here is how to compute the updated value $m^{(i+1)}$ given $\lambda = \lambda^{(i)}$ and $u = u_k$. Towards this, from \cref{eq:optProbFE3}, when the variables $\lambda$ and $u$ are fixed, the optimization problem on $m$ becomes
\begin{equation}\label{eq:optProbFE4}
	\begin{aligned}
		&\qquad m^{(i+1)} = \argmin_{m\in \calM}  && \hat{J}(m, \lambda^{(i)}, u_k) = \int_{\Omega} \left[ \frac{e_k}{m} + u_k^4 \right] \dd \bx + \lambda^{(i)} \left( \int_{\Omega} m\dd \bx - \eta |\Omega| \right)\,, \\
		&\qquad \text{and }\qquad  &&0 < m_{\text{lw}} \leq m \leq 1\,,
	\end{aligned}
\end{equation}
where $e_k = \bq(m_k) \cdot \bq(m_k) = (m_k\nabla u_k) \cdot (m_k \nabla u_k)$. Taking the variation of $\hat{J}(\cdot, \lambda^{(i)}, u_k)$ in the direction of arbitrary $w\in \calM$, and setting it to zero, gives
\begin{equation*}
	\int_{\Omega} \left[ -\frac{e_k}{^2} + \lambda^{(i)} \right] w \dd \bx = 0, \qquad \forall w\in \calM \qquad \Rightarrow \qquad m = \sqrt{\frac{e_k}{\lambda^{(i)}}}\,.
\end{equation*}
Thus, the formula for updating $m$ given $u_k$ and $\lambda^{(i)}$ is given by
\begin{equation}\label{eq:mUpdate}
	m^{(i+1)} = \min\left\{1, \max\left\{m_{\text{lw}}, \sqrt{\frac{e_k}{\lambda^{(i)}}} \right\}\right\} \,,
\end{equation}
where the upper and lower bound constraints on $m$ are enforced strongly. The algorithm for the inner iteration where $\lambda$ and $m$ are successively updated for a given outer iteration step $k$ is presented in \cref{alg:optProbInner}. This algorithm is based on the bisection method; see \cite{kumar2021direct}. 

\begin{algorithm}[H]
	\caption{Inner iteration for pair $(m_{k+1}, \lambda_{k+1})$ given $m_k, \lambda_k, u_k$.
	}
	\label{alg:optProbInner}
		\SetAlgoLined
		{\bf Input:}\\
		(a) $\eta$, lower bound $m_{\text{lw}}$, and tolerance $m_{\text{tol}}$ to check volumetric constraint \;
		(b) outer iteration step $k$ and solutions $m_k, \lambda_k, u_k$ \; 
		{\bf Setup:}\\
		(a) let $m^{(0)} = m_k$, $\lambda^{(0)} = \lambda_k$, and let $\bar{m}^{(i)}$ denote the volume average of $m^{(i)}$ \;
		(b) let $m$ and $\lambda$ are current values of variables, and let $\lambda_{min} = 0$ and $\lambda_{max} = 0$ \;
		{\bf Inner iteration (bracketing):}\\
		\If{$\bar{m}^{(0)} < \eta$}{
			$\lambda_{min} = \lambda^{(0)}$, $i=0$ \;
			\While{$\bar{m}^{(i)} < \eta$}{
				(1) update $\lambda$: $\lambda^{(i+1)} = \frac{\lambda^{(i)}}{2}$ \;
				(2) update $m$: given $m_k$, $u_k$, and $\lambda^{(i+1)}$, compute $m^{(i+1)}$ using \eqref{eq:mUpdate} \;
				(3) $m^{(i)} = m^{(i+1)}$, $\lambda^{(i)} = \lambda^{(i+1)}$, and $i \leftarrow i+1$ \;
			}
			$\lambda = 	\lambda^{(i)}$, $m = m^{(i)}$, and $\lambda_{max} = \lambda$ \;
		}
		\Else{
			$\lambda_{max} = \lambda^{(0)}$, $i=0$ \;
			\While{$\bar{m}^{(i)} > \eta$}{
				(1) update $\lambda$: $\lambda^{(i+1)} = 2\lambda^{(i)}$ \;
				(2) update $m$: given $m_k$, $u_k$, and $\lambda^{(i+1)}$, compute $m^{(i+1)}$ using \eqref{eq:mUpdate} \;
				(3) $m^{(i)} = m^{(i+1)}$, $\lambda^{(i)} = \lambda^{(i+1)}$, $i \leftarrow i+1$ \;
			}
			$\lambda = 	\lambda^{(i)}$, $m = m^{(i)}$,  and $\lambda_{min} = \lambda$ \;
		}
		{\bf Inner iteration (bisection):}\\
		$i = 0$, $\lambda^{(0)} = \lambda$, and $m^{(0)} = m$ \;
		\While{${|\bar{m}^{(0)} - \eta|} > \eta\, m_{\text{tol}} $}{
			(1) update $\lambda$: $\lambda^{(i+1)} = \frac{\lambda_{min} + \lambda_{max}}{2}$ \;
			(2) update $m$: given $m_k$, $u_k$, and $\lambda^{(i+1)}$, compute $m^{(i+1)}$ using \eqref{eq:mUpdate} \;
			(3) {update } $\lambda_{min}$ and $\lambda_{max}$ {\bf :} \\
			{
				\lIf{$\bar{m}^{(i+1)} < \eta$}
				{
					$\lambda_{min} = \lambda^{(i+1)}$
				}
				\lElse
				{
					$\lambda_{max} = \lambda^{(i+1)}$
				}
			}	
			(4) $m^{(i)} = m^{(i+1)}$, $\lambda^{(i)} = \lambda^{(i+1)}$, and $i \leftarrow i+1$ \;
		}
		$m_{k+1} = m^{(i)}$, $\lambda_{k+1} = \lambda^{(i)}$ \;
		{\bf Return: } $(m_{k+1}, \lambda_{k+1})$. 
\end{algorithm}

\end{document}